%% file: main.tex
\documentclass{article}

\input{sections/preamble}

\title{Wasserstein distances between ERGMs and Erd\H{o}s--R\'enyi models}
\author{Vilas Winstein \\ University of California, Berkeley \\ \texttt{vilas@berkeley.edu}}

\begin{document}

\maketitle

\begin{abstract}
Ferromagnetic exponential random graph models (ERGMs) are random graph models under which the presence of certain small structures (such as triangles)
is encouraged; they can be constructed by tilting an Erd\H{o}s--R\'enyi model by the exponential of a particular nonlinear Hamiltonian.
These models are mixtures of \emph{metastable wells} which each behave macroscopically like an Erd\H{o}s--R\'enyi model, exhibiting the
same laws of large numbers for subgraph counts \cite{chatterjee2013estimating}.
However, on the microscopic scale these metastable wells are very different from Erd\H{o}s--R\'enyi models, with the total variation distance between
the two measures tending to $1$ \cite{mukherjee2023statistics}.
In this article we clarify this situation by providing a sharp (up to constants) bound on the \emph{Hamming-Wasserstein} distance between the two models,
which is the average number of edges at which they differ, under the coupling which minimizes this average.
In particular, we show that this distance is $\Theta(n^{3/2})$, quantifying exactly how these models differ.

An upper bound of this form has appeared in the past \cite{reinert2019approximating}, but this was restricted to the \emph{subcritical}
(high-temperature) regime of parameters.
We extend this bound, using a new proof technique, to the \emph{supercritical} (low-temperature) regime, and prove a matching lower bound which
has only previously appeared in the subcritical regime of special cases of ERGMs satisfying a ``triangle-free'' condition \cite{ding2025second}.
To prove the lower bound in the presence of triangles, we introduce an approximation of the \emph{discrete derivative} of the Hamiltonian,
which controls the dynamical properties of the ERGM, in terms of local counts of triangles and wedges (two-stars) near an edge.
This approximation is the main technical and conceptual contribution of the article, and we expect it will be useful in a variety of other contexts as well.
Along the way, we also prove a bound on the \emph{marginal edge probability} under the ERGM via a new bootstrapping argument.
Such a bound has already appeared \cite{fang2025conditional}, but again only in the subcritical regime and using a different proof strategy.
\end{abstract}

\setcounter{tocdepth}{2}
\tableofcontents

\input{sections/intro}
\input{sections/inputs}

\input{sections/marginal}
\input{sections/hamiltonian}
\input{sections/wasserstein}

\bibliography{references}
\bibliographystyle{alpha}

\end{document}

%% file: sections/preamble.tex
\usepackage[margin=1in]{geometry}
\usepackage{amsmath, mathtools, amsfonts, amsthm, amssymb, cite, suffix, enumitem, bm, graphbox}
\mathtoolsset{showonlyrefs}
\usepackage[colorlinks]{hyperref}
\hypersetup{linkcolor=blue}
\usepackage{xcolor}

\usepackage{multirow}

\newtheorem{theorem}{Theorem}[section]
\newtheorem{lemma}[theorem]{Lemma}
\newtheorem{corollary}[theorem]{Corollary}

\newtheorem{proposition}[theorem]{Proposition}

\theoremstyle{definition}

\usepackage{mleftright}
\renewcommand{\left}{\mleft}
\renewcommand{\right}{\mright}

\renewcommand{\emptyset}{\varnothing}

\renewcommand{\P}{\mathbb{P}}
\newcommand{\E}{\mathbb{E}}
\newcommand{\R}{\mathbb{R}}
\newcommand{\N}{\mathbb{N}}

\newcommand{\cG}{\mathcal{G}}
\newcommand{\cM}{\mathcal{M}}
\newcommand{\cB}{\mathcal{B}}
\newcommand{\cN}{\mathcal{N}}
\newcommand{\cL}{\mathcal{L}}
\newcommand{\cV}{\mathcal{V}}

\newcommand{\sS}{\mathsf{S}}
\newcommand{\sT}{\mathsf{T}}

\newcommand{\cE}{\mathcal{E}}

\newcommand{\Var}{\operatorname{Var}}

\newcommand{\Cov}{\operatorname{Cov}}
\newcommand{\TV}{\mathrm{TV}}
\renewcommand{\d}{\mathbf{d}}
\newcommand{\ind}[1]{\mathbf{1}_{\{#1\}}}

\newcommand{\fr}{\frac}
\newcommand{\eps}{\varepsilon}
\newcommand{\Exp}[1]{\exp\left(#1\right)}

\WithSuffix\newcommand\ind*[1]{\mathbf{1}_{#1}}
\newcommand{\sse}{\subseteq}

\newcommand{\pball}{\cB_p^\eta}
\newcommand{\phalfball}{\cB_p^{\eta/2}}

\newcommand{\dtv}{\d_\TV}
\renewcommand{\dh}{\d_\mathrm{Ham}}

\newcommand{\dW}{\d_\mathrm{Was}}

\newcommand{\dWR}{\d_\mathrm{Was}^\R}
\newcommand{\dl}{\d_{\Lambda^*}}

\newcommand{\db}{\d_\square}
\newcommand{\dloc}{\dh^{\mathrm{loc}(v)}}
\newcommand{\dlocu}{\dh^{\mathrm{loc}(u)}}
\newcommand{\dloce}{\dh^{\mathrm{loc}(e)}}

\newcommand{\pto}{\xlongrightarrow{\P}}

\newcommand{\Ber}{\operatorname{Ber}}
\newcommand{\diam}{\operatorname{indiam}}

\newcommand{\ft}{\mathfrak{t}}
\newcommand{\sC}{\mathsf{C}}

\newcommand{\sv}{\mathsf{v}}
\newcommand{\se}{\mathsf{e}}
\newcommand{\sd}{\mathsf{d}}

\newcommand{\sN}{\mathsf{N}}
\newcommand{\hsN}{\widehat{\sN}}
\newcommand{\sH}{\mathsf{H}}

\renewcommand{\deg}{\mathsf{deg}}

\newcommand{\heH}{\widehat{n^2 \partial_e \sH}}

\newcommand{\fe}{\mathfrak{e}}
\newcommand{\fu}{\mathfrak{u}}
\newcommand{\fv}{\mathfrak{v}}

\newcommand{\st}{\mathsf{t}}
\renewcommand{\ss}{\mathsf{s}}

\renewcommand{\tilde}{\widetilde}
\newcommand{\edgeset}{\binom{[n]}{2}}

\let\temp\phi
\let\phi\varphi
\let\varphi\temp

%% file: sections/intro.tex
\section{Introduction}
\label{sec:intro}

Given a function $\sH$ of simple graphs $x$ on $n$ labeled vertices, a sample from the corresponding
\emph{exponential random graph model} (ERGM) is a random graph $\tilde{X}$ such that for any $x$ we have
\begin{equation}
\label{eq:gibbs}
    \P[\tilde{X} = x] \propto e^{n^2 \sH(x)}.
\end{equation}
Note that in this article, the function $\sH(x)$ (which we will call the \emph{Hamiltonian}) takes values of order $1$
for typical dense graphs $x$, and the scaling of $n^2$
ensures that the above measure is nontrivial (i.e.\ neither approximately uniform nor concentrated on the maximizers of $\sH$).

ERGMs have been used as models of complex networks with various clustering properties, and have found applications in fields
such as sociology and biology \cite{frank1986markov,holland1981exponential,wasserman1994social}.
Early nonrigorous analysis of these models was also contributed by statistical physicists
\cite{burda2004network,park2004solution,park2005solution}.
From a mathematical perspective, ERGMs are simple natural spin systems which incorporate higher-order
(i.e.\ not just pairwise) interactions, and they have
attracted the attention of probabilists, statisticians, and computer scientists who have investigated the model
from a variety of perspectives including
large deviations and phase transitions \cite{chatterjee2013estimating,radin2013phaseComplex,radin2013phaseERGM},
distributional limit theorems \cite{mukherjee2023statistics,fang2024normal,fang2025conditional,winstein2025quantitative},
mixing times and sampling \cite{bhamidi2008mixing,demuse2019mixing,reinert2019approximating,bresler2024metastable},
concentration of measure \cite{sambale2020logarithmic,ganguly2024sub,winstein2025concentration}, and stochastic localization
\cite{eldan2018decomposition,eldan2018exponential,mikulincer2024stochastic}.
Moreover, ERGMs are useful as a tool for studying random graphs \emph{constrained} to have a certain abnormal
number of triangles or other subgraphs \cite{radin2014asymptotics,lubetzky2015replica}.
Our definition of ERGMs yields dense graphs as samples, but some recent work has also studied
modifications which give sparse graphs instead \cite{yin2017asymptotics,cook2024typical}.
Many important works have been omitted from this list of citations, which is not meant to be comprehensive
but rather to display the breadth of the subject.

In the present work, as will be defined in Section \ref{sec:intro_definitions} below,
we consider Hamiltonians $\sH$ which can be written as finite linear combinations of subgraph densities
with nonnegative coefficients on all nonlinear terms; we call this the \emph{ferromagnetic regime}.
Ferromagnetism is a positive-correlation condition and it turns out to result in homogeneous large-scale behavior.
Specifically, under this assumption, the ERGM measure breaks down into a finite mixture of measures $\cM(n,p)$
(for $p$ ranging in some finite set), each of which is close, at the \emph{macroscopic} scale under the cut distance,
to an Erd\H{o}s--R\'enyi model $\cG(n,p)$ \cite{chatterjee2013estimating}.
This will be made precise in Section \ref{sec:intro_definitions_cut} below.
On the other hand, at the \emph{microscopic} scale the measures $\cM(n,p)$ and $\cG(n,p)$ are as far apart as possible, in the
sense that they have total variation distance tending to $1$ as $n \to \infty$ \cite{mukherjee2023statistics}, as will be
discussed in Section \ref{sec:intro_definitions_dtv} below.

So the cut distance is too coarse to distinguish $\cM(n,p)$ and $\cG(n,p)$, while the total variation distance
is too fine to compare them.
This motivates us to consider intermediate distances between these two extremes, and a natural candidate
is the \emph{Hamming distance} $\dh(x,y)$ which counts the number of edges at which the graphs $x$ and $y$ differ.
The main result of this article is that
\begin{equation}
\label{eq:main}
    \E[\dh(X,Y)] = \Theta(n^{3/2})
\end{equation}
under the optimal coupling between $X \sim \cM(n,p)$ and $Y \sim \cG(n,p)$ which minimizes this expectation.
To prove this, in particular for the lower bound, we take a perspective which treats the update
probabilities under the natural ERGM Glauber dynamics as random variables in their own right.
We then analyze these variables and derive a precise understanding of the dynamics in terms of the counts of just two
distinguished structures (namely triangles $\triangle$ and wedges $\wedge$)
present locally at each edge.

\input{sections/intro/definitions}
\input{sections/intro/results}
\input{sections/intro/outline}
\input{sections/intro/acknowledgements}

%% file: sections/intro/definitions.tex
\subsection{Problem setup and related works}
\label{sec:intro_definitions}

We work with $n$-vertex undirected graphs $x$ with no self-loops.
These graphs have vertex set $[n] = \{1,\dotsc,n\}$, and edge set $\subseteq \edgeset$, which is the collection of
two-vertex subsets of $[n]$, i.e.\ the edge set of the complete graph $K_n$.
We may thus identify $n$-vertex graphs $x$ with functions from $\edgeset$ to $\{0,1\}$, and we write $x(e)$ for the indicator that
the edge $e$ is in the graph $x$.

For any fixed graph $G$ and an $n$-vertex graph $x$, we let $\sN_G(x)$ denote the number of \emph{homomorphisms} of $G$ in $x$.
Homomorphisms are maps $\cV(G) \to [n]$ which take edges of $G$ to edges of $x$ (here $\cV(G)$ denotes the vertex set of $G$ and similarly
$\cE(G)$ denotes the edge set of $G$).
We do not require these homomorphisms to be injective, or, importantly, induced; in other words, non-edges do not have to map to non-edges.
For example, if $G = \wedge$ is a wedge graph (also called a two-star), then
\begin{equation}
    \sN_\wedge(x) = \sum_{v_1,v_2,v_3 \in [n]} x(\{v_1,v_2\}) \cdot x(\{v_2,v_3\}).
\end{equation}
Note that in the above expression, we interpret $x(\{v,v\})$ as zero since we consider graphs $x$ with no self-loops.
For any fixed graph $G$ we also define \emph{homomorphism densities} by
\begin{equation}
    \ft(G,x) \coloneqq \frac{\sN_G(x)}{n^{|\cV(G)|}},
\end{equation}
so that $\ft(G,x)$ is the probability that a uniformly random map $\cV(G) \to [n]$ is a homomorphism from $G$ to $x$.
For fixed finite graphs $G_0, G_1, \dotsc, G_K$ (where we always take $G_0$ to be a single edge), and a
parameter vector $\beta = (\beta_0, \beta_1, \dotsc, \beta_K) \in \R \times \R_{>0}^K$,
we define the ERGM Hamiltonian $\sH$ as follows:
\begin{equation}
\label{eq:hamiltonian}
    \sH(x) = \sum_{j=0}^K \beta_j \cdot \ft(G_j,x).
\end{equation}
We point out again our assumption that each $\beta_j$ (other than $\beta_0$) is positive; this is called the \emph{ferromagnetic} assumption
and it leads to somewhat homogeneous behavior at the macroscopic scale as will be discussed in Section \ref{sec:intro_definitions_cut} below.
Importantly, the ferromagnetic assumption will also allow us to employ an approximate version 
of the FKG inequality recently derived by \cite{fkg}.

The random graph $\tilde{X}$ has the distribution of the \emph{unconditioned} ERGM measure if \eqref{eq:gibbs} holds with this
choice of Hamiltonian $\sH$.
Later, we will use $X$ (with no tilde)
to denote a sample from the ERGM measure conditioned on a particular set, to be described in the next subsection.
Note also that frequently we will use $\sv$ and $\se$ to denote the number of vertices and edges of a graph $G$ under consideration,
and we reserve $\sv_j$ and $\se_j$ for the number of vertices and edges of the graphs $G_j$ in the ERGM specification.

\subsubsection{Indistinguishability of metastable wells from Erd\H{o}s--R\'enyi models in cut distance}
\label{sec:intro_definitions_cut}

It was shown in \cite{chatterjee2013estimating} via a large deviations principle that any ferromagnetic ERGM decomposes as a mixture
of measures, which we call \emph{metastable wells} or \emph{phase measures}, such that the conditional distribution within each metastable well is close
to a corresponding Erd\H{o}s--R\'enyi model, under the cut distance.

To explain this result, we give a very brief overview of graphons and the cut distance.
For a more complete understanding of this subject, see \cite{lovasz2012large}.
A \emph{graphon} is a symmetric measurable function $W : [0,1]^2 \to [0,1]$, which generalizes the notion of an adjacency matrix
of a graph.
In particular, for any graph $x$, there is a graphon $W_x$ representing $x$ defined by
\begin{equation}
    W_x = \sum_{\{u,v\} \in \edgeset} x(\{u,v\}) \cdot
    \left( \ind*{\left[\frac{u-1}{n}, \frac{u}{n}\right) \times \left[\frac{v-1}{n}, \frac{v}{n}\right)}
    + \ind*{\left[\frac{v-1}{n}, \frac{v}{n}\right) \times \left[\frac{u-1}{n}, \frac{u}{n}\right)} \right).
\end{equation}
Note that we consider graphons to be identical if they agree on a full-measure set after permuting $[0,1]$ with a measure-preserving
bijection, and so the map $x \mapsto W_x$ is non-injective.
Nonetheless, we will often interpret a graph $x$ as the graphon $W_x$ without further mention.

For any fixed $n$ we may sample a random $n$-vertex graph $Y$ from a graphon $W$ via the following procedure: first sample $n$ points
$s_v \in [0,1]$ uniformly and independently for $v \in [n]$.
Then let $Y(\{u,v\})$ be independent $\Ber(W(s_u,s_v))$ variables for each $\{u,v\} \in \edgeset$.
Using this notion of sampling, we may define subgraph densities for graphons by defining $\ft(G,W)$ to be the probability that
a $\sv$-vertex sample $Y$ from the graphon contains the $\sv$-vertex graph $G$ as a subgraph
(both graphs are defined on the same labeled vertex set $[\sv]$).
With this definition, $\ft(G, W_x) = \ft(G,x)$, and we may extend the definition of the Hamiltonian $\sH$ to be a function
of any graphon.

We may compare graphons with the \emph{cut distance}, defined as follows:
\begin{equation}
    \db(W,W') \coloneqq \inf_{\substack{\sigma : [0,1] \to [0,1] \\ \text{measure-preserving} \\ \text{bijection}}}
    \sup_{S, T \sse [0,1]} \left|
        \int_S \int_T \left( W(\sigma(s), \sigma(t)) - W'(s,t) \right) \,ds \,dt
    \right|.
\end{equation}
This formulation of the cut distance will not be used in the remainder of the article, but it should be noted that convergence
in cut distance is the same as convergence of all subgraph densities:
\begin{equation}
    \db(W_n, W) \to 0
    \qquad \Longleftrightarrow \qquad
    \ft(G,W_n) \to \ft(G,W) \quad \text{for all } G.
\end{equation}
We remark that this equivalence can be made quantitative \cite{lovasz2006limits,borgs2008convergent},
but this will not be relevant for the present article.

One of the main results of \cite{chatterjee2013estimating} is that an ERGM sample $\tilde{X}$ with distribution \eqref{eq:gibbs}
defined by the Hamiltonian \eqref{eq:hamiltonian} (under the ferromagnetic assumption) will, with high probability, be close to a
\emph{constant} graphon $W_p(s,t) = p$, where $p$ is a maximizer of the following \emph{free energy function} on $[0,1]$:
\begin{equation}
\label{eq:L}
    L_\beta(q) \coloneqq \sH(q) - \frac{1}{2} \left( q \log q + (1-q) \log (1-q) \right).
\end{equation}
Here, we write $\sH(q) = \sH(W_q)$, which is simply a polynomial in $q$ since $\ft(G,W_q) = q^{\se}$
if $G$ has $\se$ edges, and we remind the reader that the polynomial $\sH$ depends on the parameter vector $\beta$ implicitly.
Let us write $M_\beta$ for the set of global maximizers of $L_\beta$ on $[0,1]$; note that this is a finite set.
Then the combination of Theorems 3.2 and 4.1 from \cite{chatterjee2013estimating} show that for any $\eta > 0$
there are constants $C(\eta), c(\eta) > 0$ for which
\begin{equation}
\label{eq:LDP}
    \P \left[ \min_{p \in M_\beta} \db(\tilde{X}, W_p) > \eta \right] \leq C(\eta) \cdot e^{-c(\eta) \cdot n^2}.
\end{equation}
Let us now choose $\eta > 0$ small enough (depending on the ERGM specification) so that the cut-distance balls
\begin{equation}
    \pball \coloneqq \{ W \text{ graphon} : \db(W,W_p) < \eta \}
\end{equation}
are disjoint as $p$ ranges in $M_\beta$.
Further, for each $p \in M_\beta$, define $\cM(n,p)$ to be the
ERGM distribution conditioned on $\pball$.
Here we identify an $n$-vertex graph $x$ with its graphon representation $W_x$.
Note that we suppress the dependence of $\cM(n,p)$ on $\eta$ for notational convenience,
but the reader should keep this dependence in mind.
Because the sets $\pball$ are relatively stable under the natural dynamics,
we refer to the measures $\cM(n,p)$ as \emph{metastable wells} as mentioned above.

With these definitions, \eqref{eq:LDP} says that the ERGM distribution decomposes as a finite mixture
of the measures $\cM(n,p)$ for $p \in M_\beta$, with $e^{-\Omega(n^2)}$ error in total-variation distance.
Furthermore, by considering arbitrarily small values of $\eta$ in \eqref{eq:LDP}, we see that for $X \sim \cM(n,p)$ we have
\begin{equation}
    \db(X,W_p) \pto 0
\end{equation}
as $n \to \infty$.
It can also be easily checked that for $Y \sim \cG(n,p)$ (the Erd\H{o}s--R\'enyi model) we have
\begin{equation}
    \db(Y,W_p) \pto 0
\end{equation}
as well, and since $W_p$ is deterministic, it follows that
\begin{equation}
\label{eq:dboxsmal}
    \db(X,Y) \pto 0
\end{equation}
under \emph{any} coupling between $X \sim \cM(n,p)$ and $Y \sim \cG(n,p)$.
Thus the measures $\cM(n,p)$ and $\cG(n,p)$ are indistinguishable at the natural (macroscopic) scale of the cut distance.

\subsubsection{Incomparability in total variation distance}
\label{sec:intro_definitions_dtv}

On the other hand, at the (microscopic) scale of the total variation distance, the measures $\cM(n,p)$ and $\cG(n,p)$ are
as far apart as possible.
Namely, $\dtv(\cM(n,p), \cG(n,p)) \to 1$ as $n \to \infty$ for any $p \in M_\beta$.

This was first proved by \cite{mukherjee2023statistics} in the specific case of the \emph{two-star ERGM}, where in the definition
of the Hamiltonian \eqref{eq:hamiltonian}, we take $K = 1$ and $G_1 = \wedge$ to be a two-star or wedge graph.
In this case, \cite[Corollary 1.9]{mukherjee2023statistics} gives a statistical test for determining the parameters $\beta_0$ and
$\beta_1$ in the form of $\sqrt{n}$-consistent estimators $\hat{\beta}_0$ and $\hat{\beta}_1$ (note that different notation is
used in that article).
In other words, $\hat{\beta}_0(x)$ and $\hat{\beta}_1(x)$ are functions of an $n$-vertex graph $x$ such that, when
$\tilde{X}$ is sampled from the \emph{full} two-star ERGM distribution with parameters $\beta_0$ and $\beta_1$, we have
\begin{equation}
    \P \left[ \left| \hat{\beta}_0(\tilde{X}) - \beta_0 \right| > \lambda n^{-1/2} \right] \leq C_\lambda
\end{equation}
uniformly in $n$, for some constant $C_\lambda \to 0$ as $\lambda \to \infty$, and similarly for $\hat{\beta}_1(\tilde{X})$ estimating $\beta_1$.
Since Erd\H{o}s--R\'enyi models are also ERGMs with $\beta_1 = 0$, we may use these estimators to distinguish ERGMs from
Erd\H{o}s--R\'enyi models and we find that the total variation distance between the full two-star ERGM distribution and
\emph{any} Erd\H{o}s--R\'enyi model tends to $1$ as $n \to \infty$.
Since in this special case, the decomposition into metastable wells places equal weights on at most two distinct and symmetric wells,
this also shows that $\dtv(\cM(n,p), \cG(n,p)) \to 1$ as $n \to \infty$ in the case of two-star ERGMs.

Although we only consider the case of fixed parameters $\beta$ in the present work, it is worth remarking that the distinguishability
result of \cite{mukherjee2023statistics} was extended in \cite{bresler2018optimal} to consider two-star ERGM distributions with
$\beta_1$ varying with $n$.
Specifically, if $\beta_1 \gg n^{-1/2}$ then two-star ERGMs can be distinguished from Erd\H{o}s--R\'enyi models, and
if $\beta_1 \ll n^{-1/2}$ then they cannot.

As for more general ERGMs, we are not aware of any explicit statement in the literature of the fact that
$\dtv(\cM(n,p), \cG(n,p)) \to 1$, but we expect
that such a result would follow from the central limit theorems of \cite{fang2024normal,winstein2025quantitative,fkg},
which should allow for the construction of parameter estimators.
However, there may be technicalities here in the phase coexistence case (i.e.\ when $|M_\beta| > 1$), since the relative weights of metastable wells is still
an open problem beyond the two-star case, where symmetry implies that they have equal weights.
Additionally, a distributional limit theorem is not yet available in the critical case (to be defined in Section \ref{sec:intro_definitions_temp} below).
In any (non-critical) case, the fact that $\dtv(\cM(n,p),\cG(n,p)) \to 1$ follows as a direct corollary from our main result,
Theorem \ref{thm:wasserstein} below.

To match the form of \eqref{eq:dboxsmal}, we may rephrase the fact that $\dtv(\cM(n,p),\cG(n,p)) \to 1$
as the fact that
\begin{equation}
\label{eq:dtvbig}
    \ind{X \neq Y} \pto 1
\end{equation}
for any coupling between $X \sim \cM(n,p)$ and $Y \sim \cG(n,p)$.

\subsubsection{The Hamming-Wasserstein distance, and temperature regimes}
\label{sec:intro_definitions_temp}

As we have seen in \eqref{eq:dboxsmal} and \eqref{eq:dtvbig},
under any coupling of $X \sim \cM(n,p)$ and $Y \sim \cG(n,p)$,
\begin{equation}
    \db(X,Y) \pto 0 \qquad \text{and} \qquad \ind{X \neq Y} \pto 1.
\end{equation}
Both of these distances fail to capture the quantitative difference between $X$ and $Y$.
This suggests that we search for some distance which is coarser than $\ind{x \neq y}$ and finer than $\db(x,y)$.

One initial idea may be to quantify the rate of convergence of $\db(X,Y)$ to zero under some optimal coupling.
By the above results, we must have
\begin{equation}
    1 \gg \E\left[\db(X,Y)\right] \gtrsim n^{-2}
\end{equation}
under the optimal coupling minimizing this expectation, so it might be natural to consider $n^2 \db(x,y)$.
However, there is an even simpler distance which is also more interpretable and lies at roughly the same scale as $n^2 \db(x,y)$,
namely the \emph{Hamming distance} $\dh(x,y)$.
This is defined as the number of edges at which $x$ and $y$ differ, i.e.
\begin{equation}
    \dh(x,y) \coloneqq \sum_{e \in \edgeset} |x(e) - y(e)|.
\end{equation}
Of course, if $X \sim \cM(n,p)$ and $Y \sim \cG(n,p)$ are sampled independently, then $\dh(X,Y)$ is likely to be of order $n^2$, but
we can do better by coupling $X$ and $Y$ optimally.
So our main goal will be to understand
\begin{equation}
\label{eq:wasdef}
    \dW(\cM(n,p),\cG(n,p)) \coloneqq \inf_{\substack{X \sim \cM(n,p) \\ Y \sim \cG(n,p)}} \E \left[ \dh(X,Y) \right],
\end{equation}
which we call the \emph{Hamming-Wasserstein} distance, or just the Wasserstein distance for short.

The Hamming-Wasserstein distance between ERGMs and corresponding Erd\H{o}s--R\'enyi models was first considered by
\cite{reinert2019approximating} in the \emph{subcritical} or \emph{high-temperature} case, which is defined as the regime of
parameters $\beta$ for which the function $L_\beta$ defined in \eqref{eq:L} has a unique \emph{local}
maximizer $p$ and is strictly convex there, i.e.\ $L_\beta''(p) < 0$.
This is exactly the regime of parameters for which the \emph{Glauber dynamics}, where edges are resampled one-by-one with the correct
conditional probability, mixes rapidly \cite{bhamidi2008mixing}.
This rapid mixing is the source of many useful results, and can be used, as in \cite{reinert2019approximating}
to bound the Hamming-Wasserstein distance from above:
\begin{equation}
    \dW(\cM(n,p), \cG(n,p)) \lesssim n^{3/2}.
\end{equation}
Note that in the subcritical case, there is only ever \emph{one} metastable well, as metastable wells correspond to $p \in M_\beta$,
the set of global maximizers of $L_\beta$.
So instead of $\cM(n,p)$ we may just take the full ERGM distribution in this case.

Also in the subcritical regime, \cite{ganguly2024sub} proved a \emph{Gaussian concentration inequality} which may be used to derive an
alternative proof of a (slightly looser) upper bound of $O(n^{3/2} \sqrt{\log n})$ for $\dW(\cM(n,p),\cG(n,p))$.
This concentration inequality was subsequently extended to \emph{non-critical} metastable wells $\cM(n,p)$ in
the \emph{supercritical} regime by \cite{winstein2025concentration}.
The supercritical (or \emph{low-temperature} regime) is the regime of parameters $\beta$ for which there are multiple
\emph{local} maximizers of $L_\beta$, and $L_\beta$ is strictly convex at at least one of the \emph{global} maximizers.
We let $U_\beta \sse M_\beta$ denote the set of \emph{global} maximizers $p$ for which $L_\beta''(p) < 0$.
The upper bound of $O(n^{3/2} \sqrt{\log n})$ was also extended in \cite{winstein2025concentration} to any metastable well $\cM(n,p)$
for $p \in U_\beta$, except in the phase-coexistence case (where $|M_\beta| > 1$) when some graphs in the ERGM specification have cycles.
The missing part of the proof in the phase-coexistence non-forest case was a version of the FKG inequality within a metastable well, and this
was provided by \cite{fkg} which allows the bound of $O(n^{3/2} \sqrt{\log n})$ to be extended to all metastable wells $\cM(n,p)$
for $p \in U_\beta$.
The \emph{critical} case of $p \in M_\beta \setminus U_\beta$ has not been touched thus far, and in fact it is unclear if such a result
even holds in that case.

In the present article, we close the logarithmic gap by proving that $\dW(\cM(n,p),\cG(n,p)) \lesssim n^{3/2}$
for any $p \in U_\beta$ in the supercritical regime.
More importantly, we prove a matching lower bound, again for $p \in U_\beta$, in all parameter regimes, which has so far only appeared in
the subcritical regime of a special case of ERGMs where none of the graphs $G_j$ in the ERGM specification contain triangles \cite{ding2025second}.
Covering the cases where triangles are present requires new ideas; in particular, we introduce a new approximation for the dynamical behavior of the model
which is our main contribution and which we expect to be useful in a variety of other applications.
Along the way, we will also obtain an upper bound on the distance between the marginal edge probabilities
$\E[X(e)]$ and $\E[Y(e)] = p$, which is new in the supercritical regime.

%% file: sections/intro/results.tex
\subsection{Statement of results and some proof ideas}
\label{sec:intro_results}

Recall that $U_\beta$ denotes the set of $p \in [0,1]$ which are global maximizers of $L_\beta$ (defined in \eqref{eq:L})
and which satisfy $L_\beta''(p) < 0$.
Throughout this section, fix $p \in U_\beta$ and let $X \sim \cM(n,p)$, recalling the definition of the metastable well measure
$\cM(n,p)$ from Section \ref{sec:intro_definitions_cut}.
Note that $\cM(n,p)$ depends implicitly on a global choice of cut-distance ball radius $\eta > 0$, but we suppress this
dependence on $\eta$ from our notation as the results do not depend on $\eta$ once it is small enough.
Additionally, let $\cG(n,p)$ denote the Erd\H{o}s--R\'enyi measure.

Our main result is a sharp (up to constants) two-sided bound on the Hamming-Wasserstein distance between $\cM(n,p)$ and $\cG(n,p)$.
For the lower bound, we introduce the following nondegeneracy assumption:
\begin{equation}
\label{eq:nondegen}
    \textbf{Assumption: } \text{at least one graph } G_j \text{ in the ERGM specification is not a disjoint union of edges.}
\end{equation}
Recall that we have assumed $\beta_j > 0$ for $j \neq 0$ in our definition of the ERGM Hamiltonian, so any ERGM
satisfying the above assumption will have a nonzero term with a nondegenerate graph in its Hamiltonian.
We remark that this assumption is not necessary for the upper bound, but we do not expect the upper bound to be optimal
in the degenerate case, so we simply state our main result fully under this assumption.

\begin{theorem}
\label{thm:wasserstein}
Under the nondegeneracy assumption \eqref{eq:nondegen},
we have $\dW(\cM(n,p),\cG(n,p)) = \Theta(n^{3/2})$.
\end{theorem}

In the course of the proof of Theorem \ref{thm:wasserstein},
we will also derive the following bound on the marginal edge probability under $\cM(n,p)$.
This extends \cite[Proposition 1.1]{fang2025conditional} from the subcritical regime of parameters
to metastable wells in the supercritical regime, although we use a different line of reasoning.
We also mention that this result addresses \cite[Conjecture 2.1]{winstein2025concentration}
for the supercritical regime.

\begin{theorem}
\label{thm:marginal}
For $X \sim \cM(n,p)$, we have $\E \left[ X(e) \right] = p + O(n^{-1})$.
\end{theorem}

Let us briefly discuss the ideas which go into these two results, especially
Theorem \ref{thm:wasserstein}.
As a starting point, we note that the conditional expectation of an edge indicator
given the rest of the graph, i.e.
\begin{equation}
\label{eq:marginalintro}
    \E \left[ X(e) \middle| X(e') \text{ for all } e' \neq e \right],
\end{equation}
may be written explicitly as a function of the graph $X$.
For now, let us denote this variable by $\hat{p}_e(X)$, 
although we will shortly introduce more explicit notation.
Theorem \ref{thm:marginal} is the statement that
\begin{equation}
    \left| \E[\hat{p}_e(X)] - p \right| \lesssim n^{-1},
\end{equation}
while, as it turns out (see Lemma \ref{lem:optimality} below), Theorem \ref{thm:wasserstein}
is the statement that
\begin{equation}
    \E \left[ \left| \hat{p}_e(X) - p \right| \right] = \Theta(n^{-1/2}).
\end{equation}
Thus, in order to prove our results, in particular the lower bound in Theorem \ref{thm:wasserstein},
we will need a precise understanding of the distribution of the random variable $\hat{p}_e(X)$.
The main tool we develop is an approximate expression of $\hat{p}_e(X)$ solely in terms of two simple
quantities: the number of triangles in $X$ containing the edge $e$ and the number of wedges (two-stars)
in $X$ containing $e$.
We postpone a precise statement of this fact to Corollary \ref{cor:hajek} below, but let us now briefly
give a high-level idea of its proof.

First, by \eqref{eq:gibbs}, the conditional expectation \eqref{eq:marginalintro} can be written in terms of the Hamiltonian as follows:
\begin{equation}
    \hat{p}_e(X) = \frac{e^{n^2 \sH(X^{+e})}}{e^{n^2 \sH(X^{-e})} + e^{n^2 \sH(X^{+e})}}
    = \frac{e^{n^2 \partial_e \sH(X)}}{1 + e^{n^2 \partial_e \sH(X)}}
    = \phi(n^2 \partial_e \sH(X)).
\end{equation}
Here we have introduced the notation $x^{+e}$ and $x^{-e}$ to denote the graph $x$ with the edge $e$
added or removed respectively, as well as the discrete derivative
$\partial_e f(x) = f(x^{+e}) - f(x^{-e})$ and the function $\phi(s) = \frac{e^s}{1+e^s}$.
The expression $\phi(n^2 \partial_e \sH(X))$, while longer, is more expressive than $\hat{p}_e(X)$,
and this is how we will write this quantity in the future.
Note that the above equality actually does not hold if $X$ is \emph{near the boundary}
of the cut-distance ball $\pball$, as in such cases the probability of obtaining $X^{+ e}$ after resampling the edge $e$
may be zero rather than $e^{n^2 \sH(X^{+ e})}$ (and similarly for $X^{-e}$), but let us ignore this technicality for now.

Now, recalling the definition \eqref{eq:hamiltonian} of the Hamiltonian, we may write
\begin{equation}
\label{eq:peh_expansion}
    n^2 \partial_e \sH(x) = \sum_{j=0}^K \beta_j \frac{\partial_e \sN_G(x)}{n^{\sv_j-2}},
\end{equation}
where $\sv_j$ is the number of vertices in the graph $G_j$ in the ERGM specification.
Note that the quantity $\partial_e \sN_G(x)$ is the same as the number of homomorphisms of $G$ in $x^{+e}$ which
\emph{use the edge $e$}.
To be consistent with later notation, we introduce the symbol $\sN_G(x,e)$ to denote the same quantity, and view it
as a type of homomorphism count in its own right.

The heart of our argument is the approximation under $X \sim \cM(n,p)$ of the quantity $\sN_G(X,e)$, for any arbitrary graph $G$,
by a linear combination of $\sN_\triangle(X,e)$ and $\sN_\wedge(X,e)$, where $\triangle$ and $\wedge$
denote a triangle graph and a wedge graph (two-star) respectively.
The intuition for this approximation takes a dynamical route: consider sampling from $\cM(n,p)$ via Glauber dynamics $(X_t)$,
where at each step a uniformly random edge is selected and resampled according to the correct conditional
distribution (i.e.\ with probability of appearing given by \eqref{eq:marginalintro}).
We will compare how $\sN_G(X_t,e)$ evolves under these dynamics with the evolution of $\sN_\triangle(X_t,e)$
and $\sN_\wedge(X_t,e)$.
If an edge $e'$ is chosen which is not adjacent to $e$, then neither $\sN_\triangle(X_t,e)$
nor $\sN_\wedge(X_t,e)$ will change.
Correspondingly, while $\sN_G(X_t,e)$ may change in this situation, it will turn out to change by an amount which
is negligible compared to the change incurred if the chosen edge \emph{were} adjacent to $e$.
In that case, $\sN_G(X_t,e)$ may change by a lot, but we can match this change
precisely to a combination of changes in $\sN_\triangle(X_t,e)$ and in $\sN_\wedge(X_t,e)$.
The result of this analysis is presented in Proposition \ref{prop:hajek} below.

We remark that in past studies of exponential random graphs
\cite{sambale2020logarithmic,fang2024normal,winstein2025quantitative}, a \emph{linear} approximation
to $\sN_G(X)$, the overall homomorphism count, was utilized.
In other words, $\sN_G(X)$ was compared to some multiple of the total edge count.
A similar linear comparison for the number of homomorphisms using a particular \emph{vertex},
in terms of the degree of that vertex, was also useful in \cite{winstein2025quantitative}.
On the other hand, when considering $\sN_G(X,e)$, the number of homomorphisms using a particular \emph{edge},
the linear approximation would be a comparison with $\sN_\wedge(X,e)$ alone.
As shown by \cite{ding2025second}, one can get some mileage out of this approximation,
but only if all of the graphs in the ERGM specification have no triangles.

On the other hand, as our results show, when there are graphs containing triangles in the ERGM specification,
the $\sN_\triangle(X,e)$ term is necessary.
One can even do a simple back-of-the-envelope calculation in the case of Erd\H{o}s--R\'enyi models
$Y \sim \cG(n,p)$ to see that $\sN_\triangle(Y,e)$ is not well-approximated by any constant multiple of $\sN_\wedge(Y,e)$.
Thus, one of the main technical and conceptual contributions of this work is to handle the situation
where both $\sN_\triangle(X,e)$ and $\sN_\wedge(X,e)$ are relevant.
We mention also that
during the preparation of this article, another article \cite{fang2025conditional} was posted on arXiv which makes use
of a similar kind of approximation in the subcritical regime of parameters.
There, this approximation was used to prove a central limit theorem for the number of
wedges in an ERGM sample \emph{conditioned} on the number of edges, giving an example of the versatility of this type of nonlinear approximation.

Finally, we remark that there is actually an alternative proof for the lower bound in Theorem \ref{thm:wasserstein},
which relies on another result of the author \cite[Theorem 1.3]{winstein2025quantitative}.
This was pointed out by Chris Jones following a talk given by the author during the preparation of this article,
and we reproduce the argument in Section \ref{sec:wasserstein_alternative} with his permission.
Nevertheless, since the proof of \cite[Theorem 1.3]{winstein2025quantitative} is itself rather heavy, we view the original argument
for Theorem \ref{thm:wasserstein} given in the main body of this article as a more self-contained and ``bare hands'' proof.
In addition, the original proof has the benefit of motivating the aforementioned approximation which gives much sharper control over the ERGM Glauber dynamics,
and we expect this result (captured by Proposition \ref{prop:hajek} and Corollary \ref{cor:hajek} below) to be useful in a variety of other contexts.

%% file: sections/intro/outline.tex
\subsection{Outline of the article}
\label{sec:intro_outline}

We begin in Section \ref{sec:inputs} by collecting a variety of useful tools recently developed
in the literature for metastable wells in ERGMs, namely concentration and FKG inequalities.

Then in Section \ref{sec:marginal} we prove Theorem \ref{thm:marginal}, introducing a new
bootstrapping argument for the upper bound and also discussing why there is no corresponding
lower bound.
In Section \ref{sec:marginal_upperbound} we also provide a few lemmas which are used throughout the work.

Afterwards in Section \ref{sec:hamiltonian} we prove Proposition \ref{prop:hajek}
and Corollary \ref{cor:hajek}, approximating $\phi(n^2 \partial_e \sH(X))$ in terms of
$\sN_\triangle(X,e)$ and $\sN_\wedge(X,e)$ for $X \sim \cM(n,p)$.
We also prove a fluctuation lower bound on the approximating variable, which translates
to a fluctuation lower bound on $\phi(n^2 \partial_e \sH(X))$.

Finally, in Section \ref{sec:wasserstein} we wrap up with a quick proof of
Theorem \ref{thm:wasserstein}, using all of the results from the previous sections,
and we also present the alternative proof of the lower bound.

%% file: sections/intro/acknowledgements.tex
\subsection{Acknowledgements}
\label{sec:intro_acknowledgements}

I was partially supported by the NSF Graduate Research Fellowship grant DGE 2146752.
I would like to thank my advisor, Shirshendu Ganguly, for suggesting this project and for helpful discussions.
I would also like to thank Xiao Fang and Nathan Ross for helpful feedback on a draft of this article, as well as
Chris Jones for pointing out the alternative proof of the lower bound in Theorem \ref{thm:wasserstein}.

%% file: sections/inputs.tex
\section{Metastable mixing and its consequences}
\label{sec:inputs}

In this section we state various inputs which will be necessary in what follows.
Of primary importance is the concept of \emph{metastable mixing}, which is a powerful
tool that allows us to derive a variety of useful facts.
We briefly describe the relevant results on metastable mixing for ERGMs in
Section \ref{sec:inputs_mixing} below, and then in Section \ref{sec:inputs_consequences}
we state the consequences which we will use in the sequel.

\input{sections/inputs/mixing}
\input{sections/inputs/consequences}

%% file: sections/inputs/mixing.tex
\subsection{Metastable mixing}
\label{sec:inputs_mixing}

The study of the mixing time of the ERGM Glauber dynamics dates back to \cite{bhamidi2008mixing}
which proved rapid mixing in the subcritical (high-temperature) regime by showing that
the dynamics exhibits contraction in a ``good'' set $\Gamma_p^\eps$ (to be defined below)
which is also a well of attraction.
Later, \cite{bresler2024metastable} showed that $\Gamma_p^\eps$ is also \emph{large} under the metastable well measure $\cM(n,p)$,
leading to \emph{metastable mixing} or rapid mixing to stationarity from a \emph{warm start} with initial condition sampled from 
the Erd\H{o}s--R\'enyi model $\cG(n,p)$.
This was further analyzed by \cite{winstein2025concentration} which extracted a \emph{connected} subset $\Lambda^*$ of $\Gamma_p^\eps$
which is still large.
Connectedness was an important property in the use of path coupling starting from a not necessarily warm start,
and it will be crucial for our applications, in particular for deriving concentration inequalities as was done in
\cite{winstein2025concentration,winstein2025quantitative}.

To be precise and avoid confusion, in this section we use $(\tilde{X}_t^x)$ to denote the unconditioned ERGM Glauber dynamics started at $x$
(i.e.\ the Glauber dynamics with respect to the full unconditioned ERGM measure).
This only differs from the conditioned dynamics $(X_t^x)$ when it would attempt to leave the cut distance ball $\pball$ from the definition
of $\cM(n,p)$, in which case the unconditioned dynamics proceeds while the conditioned dynamics does not move.

To explain metastable mixing, let us first recall the ``good'' set $\Gamma_p^\eps$ defined in \cite{bhamidi2008mixing}.
Intuitively, this is the set where the ERGM Glauber dynamics behaves much like simply resampling edges with inclusion probability $p$.
To make this precise, for any finite graph $G$ with $\sv$ vertices and $\se$ edges, and any potential edge $e \in \edgeset$, let us define
\begin{equation}
    r_G(x, e) \coloneqq \left( \frac{N_G(x,e)}{2 \se n^{\sv-2}} \right)^{\frac{1}{\se-1}},
\end{equation}
recalling that $\sN_G(x,e)$ denotes the number of homomorphisms of $G$ in $x^{+e}$ which use the edge $e$.
To motivate this definition, note that if $Y \sim \cG(n,p)$ then $r_G(Y,e) \approx p$ with high probability.
Additionally, the update probabilities of the unconditioned ERGM Glauber dynamics are given exactly by
\begin{equation}
    \P \left[ \tilde{X}_1^x(e) = 1 \right] = \phi(n^2 \partial_e \sH(x)) = \phi \left( 2 \sum_{j=0}^K \beta_j \se_j \cdot r_{G_i}(x,e)^{\se_j-1} \right),
\end{equation}
recalling the definition of the Hamiltonian \eqref{eq:hamiltonian},
where the $G_j$ are the graphs in the ERGM specification (with $\sv_j$ vertices and $\se_j$ edges)
for $0 \leq j \leq K$.
We also remind the reader that the function $\phi(s) = \frac{e^s}{1+e^s}$ was introduced in Section \ref{sec:intro_results}.
One may easily check that each $p \in U_\beta$ (the set of strictly concave global maximizers of $L_\beta$ defined in \eqref{eq:L})
is an attracting fixed point of the map
\begin{equation}
    q \mapsto \phi(2 \sH'(q)) = \phi \left( 2 \sum_{j=0}^K \beta_j \se_j q^{\se_j-1}\right),
\end{equation}
so if all $r_G(x,e)$ are close to $p$ then $\tilde{X}_1^x$ should also have this property.
Thus, for any $\eps > 0$ we define
\begin{equation}
    \Gamma_p^\eps \coloneqq \left\{
        x \in \{0,1\}^{\edgeset} :
        \begin{array}{c}
        \left| r_G(x,e) - p \right| \leq \eps \text{ for all } e \in \edgeset \text{ and for all graphs} \\
        G \text{ with vertex count at most } \max \{ \sv_j : 0 \leq j \leq K \}
        \end{array}
    \right\}.
\end{equation}
With this definition, the following lemma holds for any $p \in U_\beta$.
In fact, it also holds for any $p$ which is an attracting fixed point of the map $q \mapsto \phi(2 \sH'(q))$,
which is equivalent to being a \emph{local} maximizer of the function $L_\beta$ defined in \eqref{eq:L},
although we will not use it in this generality.
Here we again use the notation $x^{+e}$ and $x^{-e}$ introduced in Section \ref{sec:intro_results}
for the graph $x$ with the edge added or removed respectively.
In addition, we note that the \emph{monotone coupling} between Glauber dynamics chains on $\{ 0, 1 \}^{\edgeset}$
is the coupling where all chains choose the same uniformly random edge to update and couple the subsequent Bernoulli
random variables for the edge indicators optimally.

\begin{lemma}[Lemma 18 of \cite{bhamidi2008mixing}]
\label{lem:contraction}
For all $p \in U_\beta$ there exist $\eps > 0$ and $\kappa > 0$ such that
if $x \in \Gamma_p^\eps$,
then for all $e \in \edgeset$, under the monotone coupling, we have
\begin{equation}
    \E \left[ \dh \left( \tilde{X}_1^{x^{+e}}, \tilde{X}_1^{x^{-e}} \right) \right] \leq 1 - \frac{\kappa}{n^2}.
\end{equation}
\end{lemma}

As mentioned above, it was shown in \cite{bresler2024metastable} that $\Gamma_p^\eps$ is also large under the phase measure $\cM(n,p)$
for any $p \in U_\beta$, leading to metastable mixing.
Here it is important that $p$ is not just a \emph{local} maximizer of $L_\beta$, as the methods of \cite{bresler2024metastable}
do not seem to apply to such cases, and it is an interesting question whether such a result can be extended to these
``small'' metastable wells.
In any case, the following result refines the result of \cite{bresler2024metastable}
to show that there is in fact a large subset of $\Gamma_p^\eps$
which is \emph{connected} when considered
as an induced subgraph of the nearest-neighbor Hamming cube $\{0,1\}^{\edgeset}$.
In this result, $\diam(\Lambda^*)$ denotes the diameter under the \emph{internal Hamming distance} $\dl$,
where $\dl(x,y)$ is the Hamming-length of the shortest path between $x$ and $y$
contained entirely within $\Lambda^*$.
In particular, the fact that $\diam(\Lambda^*) < \infty$ implies the connectivity property just described.

\begin{proposition}[Proposition 5.9 of \cite{winstein2025concentration}]
\label{prop:goodset}
For all small enough $\eps > 0$ there is a set $\Lambda^* \subseteq \phalfball \cap \Gamma_p^\eps$ with
$\diam(\Lambda^*) \leq 2 n^2$
and $\P[X \in \Lambda^*] \geq 1 - e^{-\Omega(n)}$ if $X \sim \cM(n,p)$, and such that if $\dh(x,\Lambda^*) \leq 1$, then
\begin{equation}
    \P\left[
        X_t^x \notin \Gamma_p^\eps \cap \phalfball \text{ for some } t \leq n^3
    \right] \leq e^{-\Omega(n)}.
\end{equation}
\end{proposition}

The importance of the fact that $\Lambda^* \subseteq \phalfball$ in addition to $\Gamma_p^\eps$ is that when the dynamics starts in $\phalfball$,
it will have no chance to attempt to leave $\pball$, at least for large enough $n$.
This means that if $\dh(x, \Lambda^*) \leq 1$ then $\tilde{X}_1^x$ and $X_1^x$ (one step of the unconditioned and conditioned dynamics respectively)
have exactly the same distribution, so we may use the contraction estimate given by Lemma \ref{lem:contraction} for the conditioned
dynamics.
Using this, we state a clean result which captures the mixing behavior of chains started in or near
the good connected set $\Lambda^*$ guaranteed by Proposition \ref{prop:goodset}.
Here as usual we let $X \sim \cM(n,p)$ and $(X_t^x)$ denotes the conditioned dynamics started at $x$.
We also slightly abuse notation and write $\dtv(A,B)$ for random variables $A$ and $B$ to denote the total-variation
distance between the \emph{distributions} of $A$ and $B$.

\begin{lemma}[Lemma 5.12 of \cite{winstein2025concentration}]
\label{lem:mixing}
There are constants $C, c > 0$ such that
\begin{equation}
    \dtv(X_{n^3}^x, X) \leq e^{-c n} 
\end{equation}
for all $x$ with $\dh(x,\Lambda^*) \leq 1$, and under the monotone coupling we have
\begin{equation}
    \E \left[ \dh(X_t^x, X_t^{x'}) \right] \leq C \left(1 - \frac{\kappa}{n^2} \right)^t
\end{equation}
for all $x \in \Lambda^*$ and $x'$ such that $\dh(x,x') = 1$, and for all $t \leq n^3$.
\end{lemma}

%% file: sections/inputs/consequences.tex
\subsection{Consequences of metastable mixing}
\label{sec:inputs_consequences}

Now let us state the consequences of metastable mixing which we will use in what follows.
These can be divided into two categories: first there are various concentration inequalities and related facts,
and second there is an approximate FKG inequality and corollaries thereof.

\subsubsection{Concentration inequalities}
\label{sec:inputs_consequences_concentration}

First we have a concentration inequality which works essentially out of the box for Lipschitz
observables of ERGM samples.
To state this, let us introduce the notion of a Lipschitz vector for a function
of $n$-vertex graphs $x$.
We say that $f$ has Lipschitz vector $\cL \in \R^{\edgeset}$ if $|f(x^{+e}) - f(x^{-e})| \leq \cL_e$
for all $e \in \edgeset$ and all $x \in \{0,1\}^{\edgeset}$.
Here as above, $x^{\pm e}$ denotes the graph $x$ with the edge $e$ added or removed.
If $f$ has Lipschitz vector $\cL$, we will also say that $f$ is $\cL$-Lipschitz.
In the following result, as usual, $X \sim \cM(n,p)$.

\begin{theorem}[Theorem 1.1 of \cite{winstein2025concentration}]
\label{thm:lipschitzconcentration}
There are constants $C,c > 0$ depending only on the ERGM specification and the choice of $p \in U_\beta$
such that, for all $\cL$-Lipschitz $f$ and all $\lambda \geq 0$ satisfying $\lambda \leq c \| \cL \|_1$,
\begin{equation}
    \P \left[ |f(X) - \E[f(X)]| > \lambda \right]
    \leq 2 \Exp{
        - c \cdot \max \left\{
            \frac{\lambda^2}{n^2 \| \cL \|_\infty^2},
            \frac{\lambda^2}{\| \cL \|_1 \| \cL \|_\infty} 
            - \Exp{\frac{C \lambda}{\| \cL \|_\infty} - c n}
        \right\}
    } + e^{- c n}.
\end{equation}
\end{theorem}

This result is easy to apply, but will not always yield the optimal bounds as the notion of Lipschitz
vector is quite coarse and often does not capture the behavior of the quantities we will 
consider at a precise enough level for our results.
As such, we will need to combine Theorem \ref{thm:lipschitzconcentration} with the following somewhat
longer but more flexible result, originally due to \cite{barbour2022long}, although we state a
version that was presented in \cite{winstein2025concentration,winstein2025quantitative}
which is a bit more adapted towards our applications.
The result holds in broader generality than just for ERGMs; to state it, we let $X$ be a random variable
in some finite state space $\Omega$.
We also use $(X_t^x)$ to denote a Markov chain with the distribution of $X$ as its stationary distribution,
started at $X_0^x = x$.
Finally, we use $P$ to denote the transition kernel of this Markov chain.

\begin{theorem}[Theorem 2.1 of \cite{barbour2022long}, Corollary 4.4 of \cite{winstein2025concentration}]
\label{thm:barbour}
Let $\Lambda \sse \Omega$ and $f : \Omega \to \R$ such that the following conditions hold for some
positive numbers $V, J, \eps, \delta, M$, a time $T \in \N$, and a point $z \in \Lambda$.
\begin{enumerate}[label=(\roman*)]
    \item
    \label{cond:variance}
    We have a bound on a proxy for the variance of $f$: for all $x \in \Lambda$,
    \begin{equation}
        \sum_{x' \in \Omega} P(x,x') \sum_{t=0}^{T-1} \left(
            \E[f(X_t^x)] - \E[f(X_t^{x'})] \right)^2 \leq V.
    \end{equation}
    \item
    \label{cond:smalljumps}
    We have a bound on the future expected difference, measured by $f$, between chains started at two adjacent initial states:
    for all $x \in \Lambda$ and $x' \in \Omega$ with $P(x,x') > 0$, and for all $t \geq 0$,
    \begin{equation}
        \left| \E[f(X_t^x)] - \E[f(X_t^{x'})] \right| \leq J.
    \end{equation}
    \item
    \label{cond:setlarge}
    The chain has a high probability of staying in $\Lambda$ for time $T$ when started from $z$:
    \begin{equation}
        \P[X_t^z \in \Lambda \text{ for } 0 \leq t < T] \geq 1 - \epsilon.
    \end{equation}

    \item
    \label{cond:mixing}
    The chain mixes rapidly with respect to $f$ when started from $z$:
    \begin{equation}
        \dtv(f(X_T^z), f(X)) \leq \delta.
    \end{equation}

    \item
    \label{cond:fbound}
    The function $f$ is bounded close to its mean: for all $x \in \Omega$,
    \begin{equation}
        |f(x) - \E[f(X)]| \leq M.
    \end{equation}
\end{enumerate}
Then for all $\lambda \geq 2 \delta M$ we have
\begin{equation}
    \P[|f(X) - \E[f(X)]| > \lambda] \leq 2 \exp \left(
        - \fr{(\lambda-2\delta M)^2}{2V + \fr{4}{3} J (\lambda - 2\delta M)}
    \right) + \epsilon + \delta.
\end{equation}
\end{theorem}

This result is, strictly speaking, not a corollary of metastable mixing.
Rather, we will apply it using the results of metastable mixing to obtain a bound.
One important tool that will be useful in this endeavor is the following lemma allowing
us to upgrade sets $\Pi$ which are merely large to sets $\Lambda$ with the properties
we need in order to apply Theorem \ref{thm:barbour} effectively.
Here, we return to the ERGM setting with $X \sim \cM(n,p)$ and $(X_t^x)$ being the conditional
ERGM Glauber dynamics, and we let $\Lambda^*$ denote the set from Proposition \ref{prop:goodset}.

\begin{lemma}[Lemma 4.3 of \cite{winstein2025quantitative}]
\label{lem:goodset_from_bigset}
Suppose that $\Pi \sse \pball$ with $\P[X \in \Pi] \geq 1 - \Exp{-\Omega(n^\xi)}$ for some $\xi \in (0,1)$.
Then there is a set $\Lambda \sse \Lambda^*$ such that
for all $x$ with $\dh(x,\Lambda) \leq 1$, we have
\begin{equation}
    \P[X_t^x \in \Pi \text{ for all } 0 \leq t \leq n^3] \geq 1 - \Exp{-\Omega(n^\xi)}.
\end{equation}
Moreover, there is some $z \in \Lambda$ for which
\begin{equation}
    \P[X_t^z \in \Lambda \text{ for all } 0 \leq t \leq n^3] \geq 1 - \Exp{-\Omega(n^\xi)}.
\end{equation}
\end{lemma}

In addition, we will need another result which allows us to control the dynamics more closely
and achieve the bounds we will need.
To state it, let us introduce the following \emph{local} Hamming distances which only measure
the difference between two graphs near a specified vertex $v \in [n]$ or edge $e \in \edgeset$:
\begin{equation}
\label{eq:dlocdef}
    \dloc(x,x') \coloneqq \sum_{e \ni v} |x(e) - x'(e)|,
    \qquad \qquad
    \dloce(x,x') \coloneqq \sum_{e' \cap e \neq \emptyset} |x(e') - x'(e')|.
\end{equation}
The following result and its subsequent corollary allow us to control these quantities.
As above, we let $(X_t^x)$ denote the conditional ERGM Glauber dynamics,
and again let $\Lambda^*$ denote the set from Proposition \ref{prop:goodset}.
Further, we let $x^{\oplus e}$ denote the graph $x$ with the status of edge $e$ flipped.

\begin{proposition}[Proposition 1.7 of \cite{winstein2025quantitative}]
\label{prop:dloc}
There is some constant $C > 0$ such that the following holds.
Let $v \in [n]$ and $e \in \edgeset$ with $v \notin e$.
Suppose that $x \in \Lambda^*$, and set $x' = x^{\oplus e}$.
Then for all $t \geq 0$, we have
\begin{equation}
    \E\left[\dloc(X_t^x,X_t^{x'})\right] \leq \fr{C}{n}.
\end{equation}
\end{proposition}

The following corollary is immediate, since if $e = \{u,v\}$ then
$\dloce(x,x') \leq \dlocu(x,x') + \dloc(x,x')$.

\begin{corollary}
\label{cor:dloce}
There is some constant $C > 0$ such that the following holds.
Let $e, e' \in \edgeset$ with $e \cap e' = \emptyset$.
Suppose that $x \in \Lambda^*$, and set $x' = x^{\oplus e'}$.
Then for all $t \geq 0$, we have
\begin{equation}
    \E\left[\dloce(X_t^x,X_t^{x'})\right] \leq \fr{C}{n}.
\end{equation}
\end{corollary}

\subsubsection{Approximate FKG inequality}
\label{sec:inputs_consequences_fkg}

Aside from concentration inequalities, the other important consequence of metastable mixing which
we will use repeatedly is the following approximate FKG inequality for $\cM(n,p)$.
Here, by $\| f \|_\infty$ we mean the maximum possible value of $f(x)$ for any $n$-vertex graph $x$.
We say that a function $f$ is coordinate-wise increasing if
$f(x^{+e}) \geq f(x^{-e})$ for all $x \in \{0,1\}^{\edgeset}$ and all $e \in \edgeset$.

\begin{theorem}[Theorem 1.2 of \cite{fkg}]
\label{thm:fkg}
Let $X \sim \cM(n,p)$.
Then there is some constant $c > 0$ such that for any coordinate-wise increasing $f$ and $g$,
\begin{equation}
    \Cov[f(X),g(X)] \geq - e^{-c n} \| f \|_\infty \| g \|_\infty.
\end{equation}
\end{theorem}

A useful corollary of the approximate FKG inequality is the following result which gives us 
control over the multilinear moments of edge indicators in ERGM samples.
This result was first proved by \cite{ganguly2024sub} in the subcritical parameter regime
using a covariance bound of \cite{newman1980normal,bulinski1998asymptotical} that relies on the FKG inequality.
Later, it was extended by \cite{winstein2025concentration,fkg} to metastable wells in the supercritical regime.

\begin{proposition}[Proposition 5.3 of \cite{fkg}]
\label{prop:multilinear}
Fix an integer $k$.
There is some constant $C>0$ such that if
$X \sim \cM(n,p)$ and $e_1, \dotsc, e_k \in \edgeset$ are distinct edges, then
\begin{equation}
    \left|
        \E \left[ \prod_{j=1}^k X(e_j) \right] - \E[X(e)]^k
    \right| \leq \frac{C}{n},
\end{equation}
where $e$ denotes an arbitrary edge in $\edgeset$.
\end{proposition}

Finally, as a consequence of Proposition \ref{prop:multilinear} and Theorem \ref{thm:lipschitzconcentration} we have the following a priori
bound on the marginal edge probability, which was also first proved in the subcritical case
by \cite{ganguly2024sub} and later extended to metastable wells in the supercritical regime by
\cite{winstein2025concentration,fkg}.

\begin{proposition}[Proposition 5.4 of \cite{fkg}]
\label{prop:apriorimarginal}
There is some constant $C > 0$ such that if $X \sim \cM(n,p)$, then
\begin{equation}
    \left| \E[X(e)] - p \right| \leq C \sqrt{\frac{\log n}{n}}
\end{equation}
for any $e \in \edgeset$.
\end{proposition}

%% file: sections/marginal.tex
\section{Marginal edge probability}
\label{sec:marginal}

In this section we prove Theorem \ref{thm:marginal}, which improves upon Proposition \ref{prop:apriorimarginal}
by showing that in fact
\begin{equation}
    \left| \E[X(e)] - p \right| \lesssim n^{-1},
\end{equation}
where as usual $X \sim \cM(n,p)$.
We expect this bound to be optimal in the sense that it cannot be improved in general, as evidenced
by simulations provided in \cite{winstein2025concentration}.
However, there is also no matching lower bound in general, even under a nondegeneracy assumption like
\eqref{eq:nondegen}; this will be explained in Section \ref{sec:marginal_nolowerbound}.

\input{sections/marginal/upperbound}
\input{sections/marginal/nolowerbound}

%% file: sections/marginal/upperbound.tex
\subsection{Upper bound via bootstrapping}
\label{sec:marginal_upperbound}

The proof of Theorem \ref{thm:marginal} below in this section can be viewed as a bootstrapping argument,
where a suboptimal upper bound (in this case given by Proposition \ref{prop:apriorimarginal})
is improved to the optimal bound.
To do this, we apply an appropriate Taylor expansion to $\E[X(e)] - p$, and by matching terms we find that
the suboptimal error bound is taken to a high power and thus subsumed, leaving only the optimal bound of order $n^{-1}$,
which ultimately traces back to the error appearing in Proposition \ref{prop:multilinear}.

We begin with a few lemmas concerning the behavior, for $X \sim \cM(n,p)$, of the important quantity $n^2 \partial_e \sH(X)$
which was discussed in Section \ref{sec:intro_results}.
These are applications of Proposition \ref{prop:multilinear} as well as
Theorem \ref{thm:lipschitzconcentration}, and they will be useful at multiple points in the proof
of Theorem \ref{thm:marginal} as well as later for Theorem \ref{thm:wasserstein}.

First, we approximate the expectation of $n^2 \partial_e \sH(X)$ by a polynomial function of the edge marginal $\E[X(e)]$.
Recall from Section \ref{sec:intro_definitions_cut}
that we interpret $\sH$ not only as a function of graphs but also as a polynomial function of $q \in [0,1]$, by
interpreting the subgraph density $\ft(G,q)$ as $q^\se$ if $G$ has $\se$ edges.

\begin{lemma}
\label{lem:multilinear_applied}
Let $X \sim \cM(n,p)$.
Then we have
\begin{equation}
    \left| \E \left[ n^2 \partial_e \sH(X) \right] - 2 \sH'(\E[X(e)]) \right| \lesssim n^{-1}.
\end{equation}
\end{lemma}

\begin{proof}[Proof of Lemma \ref{lem:multilinear_applied}]
Recalling the decomposition \eqref{eq:peh_expansion}, we have
\begin{equation}
\label{eq:peh_expansion_2}
    n^2 \partial_e \sH(x) = \sum_{j=0}^K \beta_j \frac{\sN_{G_j}(x,e)}{n^{\sv_j-2}},
\end{equation}
where $\sN_G(x,e)$ is the number of homomorphisms
of a graph $G$ in $x^{+e}$ which use the edge $e$.
Here we remind the reader that $G_j$ are the graphs in the ERGM specification, and $G_j$ has $\sv_j$ vertices and $\se_j$ edges.

For each graph $G$,
the number $\sN_G(x,e)$ can be written as a sum of indicators corresponding to the events that various maps $\cV(G) \to [n]$
which cover the edge $e$ are homomorphisms in $x^{+e}$.
If $G$ has $\sv$ vertices and $\se$ edges, then there are $2 \se n^{\sv-2}$ such indicators,
the factor of $2 \se$ arising from the choice of edge $\fe \in \cE(G)$ to map to $e$, as well as the choice of orientation,
i.e.\ whether $\fe = \{\fu,\fv\}$ is mapped to $e = \{u,v\}$ by $\fu \mapsto u$ and $\fv \mapsto v$ or by $\fu \mapsto v$
and $\fv \mapsto u$.
The factor of $n^{\sv-2}$ arises from the choice of where to map the remaining vertices of $G$.
All but $O(n^{\sv-3})$ of the maps in the sum are injective, and the indicator that an injective map is a homomorphism
is simply a product of $\se-1$ distinct edge variables (since the edge $e$ is always included by definition of $\sN_G(x,e)$).
Thus, Proposition \ref{prop:multilinear} yields
\begin{equation}
    \E \left[ \sN_G(X,e) \right] = 2 \se n^{\sv-2} \E[X(e)]^{\se-1} + O(n^{\sv-3}),
\end{equation}
where the error term comes from both the number of non-injective maps as well as the error of $O(n^{-1})$ in
Proposition \ref{prop:multilinear}.
Therefore, by \eqref{eq:peh_expansion_2} we find that
\begin{equation}
    \E \left[ n^2 \partial_e \sH(X) \right] = 2 \sum_{j=0}^K \beta_j \se_j \E[X(e)]^{\se_j-1}
    + O \left( n^{-1} \right),
\end{equation}
and the sum above is exactly $2 \sH'(\E[X(e)])$, finishing the proof.
\end{proof}

Next, we give an upper bound on the centered moments of $n^2 \partial_e \sH(X)$ using the Gaussian
concentration inequality afforded by Theorem \ref{thm:lipschitzconcentration}.

\begin{lemma}
\label{lem:concentration_moments}
Let $X \sim \cM(n,p)$.
Then for any $k \geq 1$ we have
\begin{equation}
    \E \left[ \left| n^2 \partial_e \sH(X) - \E[ n^2 \partial_e \sH(X) ] \right|^k \right] \lesssim n^{-k/2}.
\end{equation}
\end{lemma}

\begin{proof}[Proof of Lemma \ref{lem:concentration_moments}]
Let us again use the decomposition \eqref{eq:peh_expansion_2};
we will apply Theorem \ref{thm:lipschitzconcentration} to each of the terms and then combine the bounds.
First note that the $j=0$ term corresponding to a single edge already satisfies $\sN_{G_0}(x,e) \equiv 2$ which is
constant meaning this term contributes nothing to the centered moment, so it remains to apply Theorem \ref{thm:lipschitzconcentration}
to functions $\sN_G(x,e)$, where $G$ is a graph with $\sv \geq 3$ vertices and $\se$ edges.

So let us calculate a Lipschitz vector $\cL$ for such $\sN_G(x,e)$.
First, this quantity does not depend on $x(e)$, so we may take $\cL_e = 0$.
Next, for all $e'$ adjacent to $e$, there are $\lesssim n^{\sv-3}$ maps of $G$ in $x$ which
use both $e$ and $e'$, so we may take $\cL_{e'} \lesssim n^{\sv-3}$ for these $e'$.
Finally, for all $e''$ not adjacent to $e$, there are $\lesssim n^{\sv-4}$ maps using $e$ and $e''$,
so we may take $\cL_{e''} \lesssim n^{\sv-4}$ for these $e''$.
Thus the norms of the Lipschitz vector satisfy
\begin{equation}
    \| \cL \|_\infty \lesssim n^{\sv-3},
    \qquad \text{and} \qquad
    \| \cL \|_1 \lesssim n \cdot n^{\sv-3} + n^2 \cdot n^{\sv-4} \lesssim n^{\sv-2}.
\end{equation}
Thus, using the second expression in the maximum in the bound of Theorem \ref{thm:lipschitzconcentration},
we obtain constants $C, c > 0$ such that
\begin{equation}
    \P \left[ \left| \sN_G(X,e) - \E[\sN_G(X,e)] \right| > \lambda \right]
    \leq 2 \Exp{- c \frac{\lambda^2}{n^{2\sv-5}} + c e^{C \lambda n^{3-\sv} - c n}} + e^{- c n}
\end{equation}
for $\lambda \leq c n^{\sv-2}$.
If we consider only $\lambda \leq n^{\sv-2.1}$, by adjusting the constants the upper bound
simply becomes $2 \Exp{- c \lambda^2 n^{5-2\sv}}$ since the other terms are negligible.
Thus we have
\begin{align}
    \E \left[
        \left|
            \sN_G(X,e) - \E[\sN_G(X,e)]
        \right|^k
    \right] &= k \int_0^\infty t^{k-1} \P \left[ \left| \sN_G(X,e) - \E[\sN_G(X,e)] \right| > t \right] \,dt \\
    &\leq 2k \int_0^{n^{\sv-2.1}} t^{k-1} \Exp{- c t^2 n^{5-2\sv}} \,dt + 2k \int_{n^{\sv-2.1}}^{O(n^{\sv-2})} t^{k-1} \Exp{- c n^{0.8}} \,dt.
\end{align}
Here the $O(n^{\sv-2})$ in the upper bound of the second integral is any universal upper bound for $\sN_G(x,e)$
and in the second integral we have just plugged in the probability that the difference is bigger than $n^{\sv-2.1}$ instead
of bigger than $t$, giving an upper bound on the latter quantity.
Thus the second integral is simply $\Exp{-\Omega(n^{0.8})}$ as the bounds are only polynomially large.
As for the first integral, it is bounded up to constants as
\begin{equation}
    \lesssim \int_{-\infty}^\infty t^{k-1} \Exp{- c t^2 n^{5 - 2 \sv}} \,dt
    \lesssim \sqrt{n^{2 \sv - 5}} \cdot \E \left[ \left( \sqrt{n^{2\sv-5}} Z \right)^{k-1} \right],
\end{equation}
where $Z$ is a standard normal random variable.
The right-hand side above is $\lesssim n^{(\sv-2.5)k}$, and thus we find that
\begin{equation}
    \E \left[
        \left|
            \sN_G(X,e) - \E[\sN_G(X,e)]
        \right|^k
    \right] \lesssim n^{(\sv-2.5)k}.
\end{equation}
Now let us apply this to the sum \eqref{eq:peh_expansion_2}.
Since for any finite sequence of positive numbers $(a_j)_{j=1}^m$ we have
\begin{equation}
\label{eq:summoment}
    \left(\sum_{j=1}^m a_j \right)^k \leq m^k \left( \max_{1 \leq j \leq m} a_j \right)^k \leq m^k \sum_{j=1}^m a_j^k,
\end{equation}
we obtain
\begin{equation}
    \E \left[
        \left|
            n^2 \partial_e \sH(X) - \E[n^2 \partial_e \sH(X)]
        \right|^k
    \right]
    \lesssim (K+1)^k \sum_{j=0}^K \beta_j^k n^{2 k - \sv_j k} n^{(\sv_j-2.5)k} \lesssim n^{-0.5 k},
\end{equation}
which finishes the proof.
\end{proof}

Let us also combine Lemmas \ref{lem:multilinear_applied} and \ref{lem:concentration_moments} to obtain the following
bound which will be frequently used in the proof of Theorem \ref{thm:marginal}.

\begin{lemma}
\label{lem:concentration_multilinear}
Let $X \sim \cM(n,p)$.
Then for any $k \geq 1$ we have
\begin{equation}
    \E \left[ \left| n^2 \partial_e \sH(X) - 2 \sH'(\E[X(e)]) \right|^k \right] \lesssim n^{-k/2}.
\end{equation}
\end{lemma}

\begin{proof}
First note that by \eqref{eq:summoment} again we have
\begin{equation}
    \left| n^2 \partial_e \sH(X) - 2 \sH'(\E[X(e)]) \right|^k
    \leq 2^k \left| n^2 \partial_e \sH(X) - \E[n^2 \partial_e \sH(X)] \right|^k
    + 2^k \left| \E[n^2 \partial_e \sH(X)] - 2 \sH'(\E[X(e)]) \right|^k.
\end{equation}
Now Lemma \ref{lem:multilinear_applied} shows that the second term is $\lesssim n^{-k}$
and Lemma \ref{lem:concentration_moments} shows that the expectation of the first term is $\lesssim n^{-k/2}$,
giving the result.
\end{proof}

With Lemmas \ref{lem:multilinear_applied}, \ref{lem:concentration_moments}, and \ref{lem:concentration_multilinear} in hand,
we may now prove Theorem \ref{thm:marginal} using a bootstrapping argument
to improve the a priori bound of $\sqrt{n^{-1} \log n}$ to an optimal bound of $n^{-1}$.

\begin{proof}[Proof of Theorem \ref{thm:marginal}]
First, as mentioned in the discussion in Section \ref{sec:intro_results}, with $\phi(s) = \frac{e^s}{1+e^s}$, we have
\begin{equation}
\label{eq:marginalphi}
    \E[X(e)] = \E \left[ \E \left[ X(e) \middle| X(e') \text{ for } e' \neq e \right] \right]
    = \E \left[ \phi(n^2 \partial_e \sH(X)) \right] \pm e^{-\Omega(n^2)}.
\end{equation}
To verify this, note that the conditional expectation of $X(e)$ given all other edge variables is
indeed $\phi(n^2 \partial_e \sH(X))$ as long as $X \in \phalfball$, which happens with probability at least
$1 - e^{-\Omega(n^2)}$ by \eqref{eq:LDP}.
So, since $p = \phi(2\sH'(p))$, it suffices to show that
\begin{equation}
    \E \left[ \phi(n^2 \partial_e \sH(X)) - \phi(2\sH'(p)) \right] \lesssim n^{-1}.
\end{equation}
Let us begin by Taylor expanding $\phi$ around $2\sH'(p)$ to third order as follows:
\begin{align}
\label{eq:marginaltaylor1}
    \phi(n^2 \partial_e \sH(X)) - \phi(2\sH'(p))
    &= \phi'(2\sH'(p)) \cdot \left( n^2 \partial_e \sH(X) - 2\sH'(p) \right) \\
\label{eq:marginaltaylor2}
    &\qquad + \frac{\phi''(2\sH'(p))}{2} \cdot \left( n^2 \partial_e \sH(X) - 2\sH'(p) \right)^2 \\
\label{eq:marginaltaylor3}
    &\qquad + O \left( \left| n^2 \partial_e \sH(X) - 2\sH'(p) \right|^3 \right).
\end{align}
We will take the expectation of the above quantity
and then replace all instances of $n^2 \partial_e \sH(X)$ by $2\sH'(\E[X(e)])$ and track the
errors which accumulate.
Our goal with this is to show that $\E[X(e)]$ is approximately a fixed point of the map $q \mapsto \phi(2 \sH'(q))$,
but that map has a unique fixed point in the window we are considering, namely $p$.
This will allow us to conclude that $\E[X(e)]$ is close to $p$.

Let us first expand the second factor in \eqref{eq:marginaltaylor1} as
\begin{align}
\label{eq:marginal11}
    n^2 \partial_e \sH(X) - 2\sH'(p)
    &= \left( n^2 \partial_e \sH(X) - 2\sH'(\E[X(e)]) \right) \\
\label{eq:marginal12}
    &\qquad + \left( 2\sH'(\E[X(e)]) - 2\sH'(p) \right),
\end{align}
and note that the expectation of \eqref{eq:marginal11} is $\lesssim n^{-1}$ by
Lemma \ref{lem:multilinear_applied}.
So since $\phi'(s) \leq \frac{1}{4}$ for all $s \in \R$, we obtain
\begin{equation}
\label{eq:marginal1end}
    \E \left[ \eqref{eq:marginaltaylor1} \right] = \phi'(2\sH'(p)) \cdot \left( 2\sH'(\E[X(e)]) - 2\sH'(p) \right)
    + O \left( n^{-1} \right).
\end{equation}

Similarly, we expand the second factor in \eqref{eq:marginaltaylor2} as
\begin{align}
\label{eq:marginal21}
    \left( n^2 \partial_e \sH(X) - 2\sH'(p) \right)^2
    &= \left(n^2 \partial_e \sH(X) - 2\sH'(\E[X(e)]) \right)^2 \\
\label{eq:marginal22}
    &\qquad + 2 \left(n^2 \partial_e \sH(X) - 2\sH'(\E[X(e)]) \right) \cdot \left( 2\sH'(\E[X(e)]) - 2\sH'(p) \right) \\
\label{eq:marginal23}
    &\qquad + \left( 2\sH'(\E[X(e)]) - 2\sH'(p) \right)^2.
\end{align}
Now the expectation of \eqref{eq:marginal21} is $\lesssim n^{-1}$ by Lemma \ref{lem:concentration_multilinear},
and the expectation of \eqref{eq:marginal22} is also $\lesssim n^{-1}$ by Lemma \ref{lem:multilinear_applied} as in the previous paragraph,
noting that the second factor of \eqref{eq:marginal22} is deterministic and of order at most $\sqrt{n^{-1} \log n} \lesssim 1$ by our a priori bound.
So, since $\phi''(s) \lesssim 1$ for $s \in \R$, we obtain
\begin{equation}
\label{eq:marginal2end}
    \E \left[ \eqref{eq:marginaltaylor2} \right]
    = \frac{\phi''(2\sH'(p))}{2} \cdot \left( 2\sH'(\E[X(e)]) - 2 \sH'(p) \right)^2 + O \left( n^{-1} \right).
\end{equation}

Finally, using the inequality $|a+b|^3 \leq 8|a|^3 + 8 |b|^3$ as in \eqref{eq:summoment}, we bound \eqref{eq:marginaltaylor3} as 
\begin{align}
\label{eq:marginal31}
    \left| n^2 \partial_e \sH(X) - 2 \sH'(p) \right|^3
    &\leq 8 \left|n^2 \partial_e \sH(X) - 2\sH'(\E[X(e)]) \right|^3 \\
\label{eq:marginal32}
    &\qquad + 8 \left| 2\sH'(\E[X(e)]) - 2\sH'(p) \right|^3.
\end{align}
The expectation of \eqref{eq:marginal31} is $\lesssim n^{-3/2}$ by Lemma \ref{lem:concentration_multilinear}.
And since $\sH'$ is a polynomial it is Lipschitz on the interval $[0,1]$, so by our a priori marginal
bound in Proposition \ref{prop:apriorimarginal} which states that $\E[X(e)] = p + O(\sqrt{n^{-1} \log n})$,
we find that \eqref{eq:marginal32} is $\lesssim n^{-3/2} \log^{3/2}(n)$.
Thus, overall we have
\begin{equation}
\label{eq:marginal3end}
    \E \left[ \eqref{eq:marginaltaylor3} \right] \lesssim n^{-3/2} \log^{3/2}(n).
\end{equation}

Combining \eqref{eq:marginal1end}, \eqref{eq:marginal2end}, and \eqref{eq:marginal3end}, we obtain
\begin{align}
    \E\left[\phi(n^2 \partial_e \sH(X)) - \phi(2\sH'(p))\right]
    &= \phi'(2\sH'(p)) \cdot \left( 2 \sH'(\E[X(e)]) - 2\sH'(p) \right) \\
    &\qquad + \frac{\phi''(2\sH'(p))}{2} \cdot \left( 2 \sH'(\E[X(e)]) - 2\sH'(p) \right)^2 
    + O \left( n^{-1} \right).
\end{align}
Now by another Taylor expansion, the right-hand side above is
\begin{equation}
    \phi(2\sH'(\E[X(e)])) - \phi(2\sH'(p)) + O \left( \left| 2 \sH'(\E[X(e)]) - 2 \sH'(p) \right|^3 \right) + O(n^{-1}),
\end{equation}
and the error term in the middle is also $n^{-3/2} \log^{3/2}(n)$ by the a priori bound of
Proposition \ref{prop:apriorimarginal} and the fact that $\sH'$ is Lipschitz.
So we find that
\begin{equation}
    \E\left[\phi(n^2 \partial_e \sH(X)) - \phi(2\sH'(p))\right]
    = 
    \phi(2\sH'(\E[X(e)])) - \phi(2\sH'(p)) + O(n^{-1}),
\end{equation}
and subtracting $\phi(2\sH'(p))$ from both sides yields
\begin{equation}
    \E\left[\phi(n^2 \partial_e \sH(X))\right] = \phi(2\sH'(\E[X(e)])) + O \left( n^{-1} \right).
\end{equation}
Since the left-hand side above is $\E[X(e)] \pm e^{-\Omega(n^2)}$ by \eqref{eq:marginalphi},
we find that
\begin{equation}
\label{eq:marginalfixedpoint}
    \left| \E[X(e)] - \phi(2\sH'(\E[X(e)])) \right| \lesssim n^{-1}.
\end{equation}
Now, the map $q \mapsto q - \phi(2\sH'(q))$ has a unique zero, namely $p$, in the interval
$[ p - \eta, p + \eta]$ for $\eta$ small enough, and for $X \sim \cM(n,p)$ we have,
a priori, that $\E[X(e)] \in [ p - \eta, p + \eta]$ (recall the definition of $\cM(n,p)$ involved
a cut-distance error $\eta > 0$).
Moreover, the derivative of this map is strictly positive on this interval as well, for small enough $\eta>0$,
since $\frac{d}{dq} \phi(2\sH'(q)) < 1$ at $q=p$.
Thus \eqref{eq:marginalfixedpoint} implies that
\begin{equation}
    \left| \E[X(e)] - p \right| \lesssim n^{-1},
\end{equation}
finishing the proof.
\end{proof}

With Theorem \ref{thm:marginal} in place, we present a corollary
akin to the lemmas earlier in the section
which will be useful at various points
throughout the proof of Theorem \ref{thm:wasserstein}.

\begin{corollary}
\label{cor:concentration_marginal}
Let $X \sim \cM(n,p)$.
Then for any $k \geq 1$ we have
\begin{equation}
    \E \left[ \left| n^2 \partial_e \sH(X) - 2 \sH'(p) \right|^k \right] \lesssim n^{-k/2}.
\end{equation}
\end{corollary}

\begin{proof}[Proof of Corollary \ref{cor:concentration_marginal}]
Let us expand the argument as
\begin{equation}
    \left| n^2 \partial_e \sH(X) - 2 \sH'(p) \right|^k
    \leq 2^k \left| n^2 \partial_e \sH(X) - 2 \sH'(\E[X(e)]) \right|^k
    + 2^k \left| 2 \sH'(\E[X(e)]) - 2 \sH'(p) \right|^k.
\end{equation}
First, by Lemma \ref{lem:concentration_multilinear} the expectation of the first term above is $\lesssim n^{-k/2}$.
Second, since $\sH'$ is Lipschitz, Theorem \ref{thm:marginal} implies that the second term is $\lesssim n^{-k} \leq n^{-k/2}$,
finishing the proof.
\end{proof}

%% file: sections/marginal/nolowerbound.tex
\subsection{There is no lower bound}
\label{sec:marginal_nolowerbound}

In this section we explain why there is no lower bound for $|\E[X(e)] - p|$ matching the upper bound
given in Theorem \ref{thm:marginal}.
In short, \cite{fang2025conditional} showed that in the subcritical regime of parameters
there is an explicit formula for a constant $c^*$ such that
\begin{equation}
\label{eq:marginalformula}
    \E[X(e)] = p + c^* n^{-1} + O\left( n^{-3/2} \right).
\end{equation}
We will use this explicit formula to prove that $c^*$ may be zero.

In order to state the explicit formula for $c^*$, recall that
$G_0, \dotsc, G_K$ be the graphs in the ERGM specification.
Let us also use $\sv_l, \se_l, \ss_l, \st_l$ to denote the number of vertices, edges,
wedges (i.e.\ two-stars) and triangles in $G_l$ respectively.
Additionally, we fix an ordering $\fe_l^1, \dotsc \fe_l^{\se_l}$ of the edges in $G_l$
and let $\ss_l^j$ and $\st_l^j$ denote the number of two-stars and triangles in $G_l \setminus \fe_l^j$.
Also set
\begin{equation}
    \sS = \sum_{l=1}^K \beta_l \ss_l p^{\se_l}
    \qquad \text{and} \qquad
    \sT = \sum_{l=1}^K \beta_l \st_l p^{\se_l}
\end{equation}
Then the following value of $c^*$ satisfies \eqref{eq:marginalformula} in the subcritical regime
by \cite[Lemma B.1]{fang2025conditional}:
\begin{align}
    c^* = \fr{p(1-p)}{1 - \phi'(p)} \Bigg[
        &\sum_{l=1}^K \beta_l \sum_{j=1}^{\se_l} \left(
            4 \ss_l^j p^{\se_l - 3}
            \fr{
                (1-p)^2 \sS
            }{
                1 - 2 p^{-1} (1-p) \sS
            }
            + 12 \st_l^j p^{\se_l-4}(1-p)^3 \sS
        \right) \\
        &\qquad + (1-2p) \left(
            \fr{
                8 p^{-3} (1-p) S^2
            }{
                1 - 2 p^{-1} (1-p) \sS
            }
            + 36 p^{-4} (1-p)^2 \sT^2
        \right) \\
        &\qquad \qquad - \sum_{l=1}^K \beta_l \se_l p^{\se_l - 1} (\sv_l - 2) (\sv_l - 3)
    \Bigg].
\end{align}
Note that when $\beta_0, \dotsc, \beta_K$ are in the subcritical regime of parameters,
the optimal density $p$ is an analytic function of the parameters; thus, in the subcritical
regime, the value $c^*$ is also.
A stronger condition than being subcritical is the \emph{Dobrushin regime} defined by
\begin{equation}
    \sH''(1) < 2.
\end{equation}
This condition implies subcriticality and thus analyticity of $c^*$ as a function of the parameters
\cite{bianchi2024limit}.

Now, in the case where $G_0$ is an edge, $G_1$ is a two-star and $G_2$ is a triangle, 
$\ss_1 = 1$, $\st_1 = 0$, $\ss_2 = 3$ and $\st_2 = 1$.
Additionally, $\ss_1^j = \st_1^j = 0$ for $j=1,2$ and $\ss_2^j = 1$ and $\st_2^j = 0$ for $j=1,2$.
Thus in this case we have
\begin{equation}
    \sS = \beta_1 p^2 + 3 \beta_2 p^3
    \qquad \text{and} \qquad
    \sT = \beta_2 p^3,
\end{equation}
and so we have
\begin{align}
    c^* = \fr{p(1-p)}{1-\phi'(p)} \Bigg[
        &\beta_1 \cdot 3 \cdot 4 \cdot 1 \cdot p^{-1} \fr{(1-p)^2 (\beta_1 p^2 + 3 \beta_2 p^3)}{1 - 2 p^{-1}(1-p) (\beta_1 p^2 + 3 \beta_2 p^3)} \\
        &+ (1 - 2p) \left(
            \fr{8 p^{-3} (1-p) (\beta_1 p^2 + 3 \beta_2 p^3)^2}
            {1 - 2p^{-1} (1-p) (\beta_1 p^2 + 3 \beta_2 p^3)} 
            + 36 p^{-4} (1-p)^2 \beta_2^2 p^6
        \right)
    \Bigg].
\end{align}
Note that the final term disappears since $\sv_1 = \sv_2 = 3$.

Now, as can be checked numerically or otherwise, writing $c^* = c^*(\beta_0,\beta_1,\beta_2)$, we have
\begin{equation}
    c^*(0,0.01,0.01) \approx + 0.000069
    \qquad \text{and} \qquad
    c^*(0,0.1,0.1) \approx - 0.119409.
\end{equation}
Thus, as long as $c^*$ is a continuous function of the parameters between these two values,
there must be some parameter values for which $c^* = 0$.
Let us check that these parameter values in fact all lie in the Dobrushin regime.
We have
\begin{equation}
    \sH(p) = \beta_0 p + \beta_1 p^2 + \beta_2 p^3,
    \qquad \text{and so} \qquad
    \sH''(1) = 2 \beta_1 + 6 \beta_2.
\end{equation}
Thus for all $\beta_1, \beta_2 \in [0.01, 0.1]$, we do indeed have $\sH''(1) < 2$,
meaning we are in the Dobrushin regime and $c^*$ is an analytic function of the parameters.

Note also that one could have examined the simulations provided in \cite[Section 2]{winstein2025concentration},
and there found that it is possible for $c^*$ to be both positive and negative, at least experimentally,
for different parameter values $\beta$ with the same graphs in the ERGM specification.
Following a path through the phase uniqueness regime of the parameter space which connects the parameters
giving a positive value of $c^*$ with the parameters giving a negative value would yield the same conclusion.

%% file: sections/hamiltonian.tex
\section{Fluctuations of the Hamiltonian differential}
\label{sec:hamiltonian}

In this section we come to the main technical and conceptual contribution of the article, a
lower bound of order $n^{-1/2}$ on the fluctuations of $n^2 \partial_e \sH(X)$ for $X \sim \cM(n,p)$,
which is critical for our proof of the lower bound in Theorem \ref{thm:wasserstein}.
The proof will proceed via first reducing $n^2 \partial_e \sH(X)$ to a simpler (but still nonlinear) function of $X$
which depends only on the very local structure near $e$, namely the number of wedges and triangles adjacent to a particular edge.
This will be accomplished in Section \ref{sec:hamiltonian_hajek}.
Then in Section \ref{sec:hamiltonian_fkg} we will give a lower bound on the fluctuations of this simpler function, and translate
that to a lower bound on the fluctuations of $n^2 \partial_e \sH(X)$.

We note that during the preparation of this article, a similar comparison of $n^2 \partial_e \sH(X)$
with the number of wedges and triangles adjacent to an edge appeared in \cite{fang2025conditional}
in the subcritical regime, and this was used for the purpose of proving a central limit theorem for the total number
of wedges in $X$, conditioned on the total number of edges in $X$.
However, our proof techniques in the present work are rather different and, beyond working for the supercritical regime,
they give improved tail bounds even in the subcritical regime.
This is due to the strategy employed, which is at a high level similar in spirit to a method introuced in
\cite{winstein2025quantitative}
that improves upon tail bounds obtainable by the higher-order concentration methods given in \cite{sambale2020logarithmic}
which are used in \cite{fang2024normal,fang2025conditional} in the subcritical regime.
For more information on the comparison between these two methods,
the reader may consult \cite[Section 1.3.3]{winstein2025quantitative}.

\input{sections/hamiltonian/hajek}
\input{sections/hamiltonian/fkg}

%% file: sections/hamiltonian/hajek.tex
\subsection{Reduction of Hamiltonian differential to a simpler variable}
\label{sec:hamiltonian_hajek}

We begin by recalling from \eqref{eq:peh_expansion} and \eqref{eq:peh_expansion_2} the decomposition
\begin{equation}
\label{eq:peh_expansion_3}
    n^2 \partial_e \sH(x) = \sum_{j=0}^K \beta_j \frac{\sN_G(x,e)}{n^{\sv_j-2}},
\end{equation}
where again we use $\sN_G(x,e)$ to denote the number of homomorphisms of $G$ in $x^{+e}$ which use the edge $e$; this is equal to
$\partial_e \sN_G(x)$.
The main result of this section is the following proposition, which reduces the understanding of the fluctuations
(for $X \sim \cM(n,p)$) of
$\sN_G(X,e)$ to those only of $\sN_\triangle(X,e)$ and $\sN_\wedge(X,e)$, where $\triangle$ and $\wedge$ represent
a triangle graph and a wedge graph respectively.
Later we will use this result with \eqref{eq:peh_expansion_3} to derive the fluctuation lower bound on $n^2 \partial_e \sH(X)$
in Corollary \ref{cor:hajek}.

\begin{proposition}
\label{prop:hajek}
For any graph $G$ with $\sv$ vertices and $\se$ edges,
there are constants $\sC_G^\triangle, \sC_G^\wedge \geq 0$, with at least one of them
strictly positive if $G$ is not a disjoint union of edges, such that the following holds.
Define
\begin{equation}
    \hsN_G(x,e) \coloneqq \sN_G(x,e)
    - \sC_G^\triangle \cdot p^{\se-3} n^{\sv-3} \cdot \sN_\triangle(x,e)
    - \sC_G^\wedge \cdot p^{\se-2} n^{\sv-3} \cdot \sN_\wedge(x,e).
\end{equation}
Then for any $\zeta \in (0, \frac{3}{4})$ and any $\gamma > 0$ there exists some $c > 0$ such that if $X \sim \cM(n,p)$ then
\begin{equation}
    \P \left[
        \left|
            \hsN_G(X,e) - \E[\hsN_G(X,e)]
        \right| > \gamma n^{\sv - 2.75 + \zeta}
    \right] \leq \Exp{-c n^{\frac{4}{3} \zeta}}.
\end{equation}
\end{proposition}

Let us now set up the preliminary notation and facts which we will need for the proof of Proposition \ref{prop:hajek}.
For any edge $\fe \in \cE(G)$ (the edge set of $G$),
let us define $\sN_G(x, \fe \mapsto e)$ as the number of homomorphisms of $G$ in $x^{+e}$ which map 
$\fe$ to $e$.
Note that if $\fe = \{\fu,\fv\}$ and $e = \{u,v\}$ then this is the number of homomorphisms which either
map $\fu \mapsto u$ and $\fv \mapsto v$ or map $\fu \mapsto v$ and $\fv \mapsto u$.
With this definition, we have
\begin{equation}
\label{eq:Nxedecomp}
    \sN_G(x,e) = \sum_{\fe \in \cE(G)} \sN_G(x, \fe \mapsto e)
    + O(n^{\sv-3}),
\end{equation}
since there must be some $\fe \in \cE(G)$ which maps to $e$,
and if there are multiple $\fe \in \cE(G)$ which map to $e$
then there are at most $\sv-3$ vertices of $\cG$ remaining to map to $[n]$.
In other words, there are at most $O(n^{\sv-3})$ noninjective maps from $\cV(G)$ to $[n]$ which cover $e$.

Thus, to prove Proposition \ref{prop:hajek}, it suffices to show the existence of constants
$\sC_G^\Delta(\fe), \sC_G^\wedge(\fe) \geq 0$ for each $\fe \in \cE(G)$,
with at least one of them strictly positive if $\fe$ has any adjacent edges in $G$,
such that if we define
\begin{equation}
\label{eq:hsngxeee_def}
    \hsN_G(x,\fe \mapsto e) = \sN_G(x,\fe \mapsto e)
    - \sC_G^\triangle(\fe) \cdot p^{\se-3} n^{\sv-3} \cdot \sN_\triangle(x,e)
    - \sC_G^\wedge(\fe) \cdot p^{\se-2} n^{\sv-3} \cdot \sN_\wedge(x,e)
\end{equation}
then for $X \sim \cM(n,p)$, for any $\zeta \in (0,\frac{3}{4})$, and for any $\gamma > 0$, there is a constant $c > 0$ such that
\begin{equation}
\label{eq:hajek_fe}
    \P\left[
        \left|
            \hsN_G(X,\fe \mapsto e) - \E[\hsN_G(X,\fe \mapsto e)]
        \right| > \gamma n^{\sv - 2.75 + \zeta}
    \right] < \Exp{- c n^{\frac{4}{3} \zeta}}.
\end{equation}
Indeed, the conclusion of Proposition \ref{prop:hajek} will follow from \eqref{eq:hajek_fe} by a union bound
over the finitely many edges $\fe \in \cE(G)$, expressing $\hsN_G(X,e)$ as a sum using \eqref{eq:Nxedecomp}
and adjusting the constant $c$.
To see the details of a very similar argument spelled out, the reader may consult the proof of Corollary \ref{cor:dhfluct} below.

To prove \eqref{eq:hajek_fe} we will apply a strategy which is somewhat similar to \cite[Propositions 4.1 and 4.4]{winstein2025quantitative}.
Namely, we will use Theorem \ref{thm:barbour}, meaning that we need to control the behavior of $\hsN_G(X_t,\fe \mapsto e)$
under the conditioned ERGM Glauber dynamics $(X_t)$.
We will do this in stages, first by using Theorem \ref{thm:lipschitzconcentration} to control the \emph{changes} under
the dynamics in terms of some simpler quantities involving only triangles $\triangle$ and wedges $\wedge$.
Then we will use this as an input to Theorem \ref{thm:barbour} to derive \eqref{eq:hajek_fe}.
One key difference between \cite[Propositions 4.1 and 4.4]{winstein2025quantitative} and the present work is that a \emph{linear}
approximation (i.e.\ just using $\sN_\wedge(x,\fe \mapsto e)$ and not $\sN_\triangle(x,\fe \mapsto e)$) will not suffice,
as discussed in Section \ref{sec:intro_results}.

In any case, we first we need to understand the changes of $\hsN_G(X_t, \fe \mapsto e)$ under Glauber dynamics,
which will reduce to considering the quantities $\partial_{e'} \hsN_G(X,\fe \mapsto e)$ for a typical sample $X \sim \cM(n,p)$
and for $e' \neq e$ (note that $\hsN_G(x, \fe \mapsto e)$ does not depend on $x(e)$ since in the calculation $e$ is always
added to $x$).
Let us start by considering $\partial_{e'} \sN_G(x, \fe \mapsto e)$ for general deterministic graphs $x$.
This is the number of homomorphisms of $G$ in $x^{+e+e'}$
taking $\fe$ to $e$ and also mapping some other edge in $G$ to $e'$; let us denote this number by $\sN_G(x, \fe \mapsto e, e')$.
Then by similar reasoning as to \eqref{eq:Nxedecomp}, we have
\begin{equation}
\label{eq:Nxeeedecomp}
    \sN_G(x,\fe \mapsto e, e') = \sum_{\fe' \in \cE(G) \setminus \{ \fe \} } \sN_G(x,\fe \mapsto e, \fe' \mapsto e')
    + O(n^{\sv-4}).
\end{equation}
Here $\sN_G(x, \fe \mapsto e, \fe' \mapsto e')$ counts the number of homomorphisms which map
$\fe$ to $e$ and map $\fe'$ to $e'$; again, this includes both orientations for each edge as with $\sN_G(x, \fe \mapsto e)$,
but when $\fe$ and $\fe'$ share a vertex the number of possibilities will be reduced, as the orientation of the mapping of $\fe$
will determine the orientation of the mapping of $\fe'$ (this will be relevant shortly).
The equation \eqref{eq:Nxeeedecomp} holds because there are at most $O(n^{\sv-4})$ noninjective maps counted by
$\sN_G(X,\fe \mapsto e, e')$, and this can be improved to $O(n^{\sv-5})$ if $e$ and $e'$ don't share a vertex.
So we remark that the error term can be improved to $O(n^{\sv-5})$ if $e$ and $e'$ do not share a vertex, although
this will be not important in the sequel.
The reason for this is that we will consider the case of $e \cap e' = \emptyset$ to be negligible as there are 
at most $O(n^{\sv-4})$ maps in this case anyway.
Indeed, we actually have the following bounds:
\begin{align}
    \label{eq:neighbors}
    \sN_G(x, \fe \mapsto e, e') &= \sum_{\substack{\fe' \in \cE(G) \setminus \{\fe\} \\ \fe \cap \fe' \neq \emptyset}}
    \sN_G(x,\fe \mapsto e, \fe' \mapsto e') + O(n^{\sv-4}) & & \text{if } e \cap e' \neq \emptyset, \\
    \label{eq:notneighbors}
    \sN_G(x, \fe \mapsto e, e') &= O(n^{\sv-4}) & & \text{if } e \cap e' = \emptyset.
\end{align}
So we will only need to analyze $\sN_G(x,\fe \mapsto e, \fe' \mapsto e')$ for in the case where both $\fe \cap \fe' \neq \emptyset$
and $e \cap e' \neq \emptyset$.

Now we need to consider two cases: either $\fe$ and $\fe'$ span a triangle in $G$, or not.
If they do span a triangle, then in order for a map taking $\fe \mapsto e$ and $\fe' \mapsto e'$
to be a homomorphism in $x$, the third edge $e''$ of the triangle spanned by $e,e'$ in the complete graph $K_n$
must also be present in $x$, meaning that
\begin{equation}
    \sN_G(x, \fe \mapsto e, \fe' \mapsto e') = x(e'') \cdot \sN_G(x^{+e''}, \fe \mapsto e, \fe' \mapsto e').
\end{equation}
Note that on the right-hand side the edges $e$ and $e'$ are also implicitly added to $x$ when counting
the homomorphisms, by the definition of $\sN_G(x, \fe \mapsto e, \fe' \mapsto e')$.
This leads to the following representation in the case where $e \cap e' \neq \emptyset$,
where we denote by $e''$ the edge completing the triangle with $e$ and $e'$:
\begin{equation}
\label{eq:Ngxeeerepresentation}
    \sN_G(x, \fe \mapsto e, e')
    = \sum_{\fe' \in \cE_\triangle(G,\fe)} x(e'') \cdot \sN_G(x^{+e''}, \fe \mapsto e, \fe' \mapsto e')
    + \sum_{\fe' \in \cE_\wedge(G,\fe)} \sN_G(x, \fe \mapsto e, \fe' \mapsto e') + O(n^{\sv-4}).
\end{equation}
In the above expression, $\cE_\triangle(G,\fe)$ denotes the set of edges $\fe' \in \cE(G) \setminus \{\fe\}$,
which span a triangle with $\fe$ in $G$, and $\cE_\wedge(G,\fe)$ denotes the set of edges $\fe' \in \cE(G) \setminus \{\fe\}$
which neighbor but \emph{do not} span a triangle with $\fe$ in $G$.

The first part of our analysis will be to show that the expressions in the sums in \eqref{eq:Ngxeeerepresentation}
are concentrated around deterministic values.
The following two lemmas showing this for the two types of expressions have essentially the same proof, which we provide
below.
As usual, in the following statements we let $X \sim \cM(n,p)$, and $G$ is a graph with $\sv$ vertices and $\se$ edges.

\begin{lemma}
\label{lem:stepbound_triangle}
Suppose $\fe \in \cE(G)$ and $\fe' \in \cE_\triangle(G,\fe)$, i.e.\ $\fe$ and $\fe'$ span a triangle in $G$.
For any $\alpha \in (0,\fr{1}{2})$ and $\gamma > 0$, there exists a constant $c$ such that for all $e,e' \in \edgeset$ which share a vertex
(and thus span a triangle in $K_n$ which is completed by an edge which we denote by $e''$), we have
\begin{equation}
    \P \left[
        \left|
            \sN_G(X^{+e''}, \fe \mapsto e, \fe' \mapsto e') - 2 p^{\se-3} n^{\sv-3}
        \right| > \gamma n^{\sv-3.5+\alpha}
    \right] \leq \Exp{-c n^{2\alpha}}.
\end{equation}
\end{lemma}

\begin{lemma}
\label{lem:stepbound_wedge}
Suppose $\fe \in \cE(G)$ and $\fe' \in \cE_\wedge(G,\fe)$, i.e.\ $\fe$ and $\fe'$ share a vertex but \emph{do not}
span a triangle in $G$.
For any $\alpha \in (0,\fr{1}{2})$ and $\gamma > 0$, there exists a constant $c$ such that for all $e,e' \in \edgeset$ which share a vertex,
we have
\begin{equation}
    \P \left[
        \left|
            \sN_G(X, \fe \mapsto e, \fe' \mapsto e') - 2 p^{\se-2} n^{\sv-3}
        \right| > \gamma n^{\sv-3.5+\alpha}
    \right] \leq \Exp{-c n^{2\alpha}}.
\end{equation}
\end{lemma}

Before proving these lemmas, let us see how the proof of Proposition \ref{prop:hajek} will proceed at a heuristic level.
By the lemmas and \eqref{eq:Ngxeeerepresentation}, when we change the status in $X$ of an edge $e'$ adjacent (but not equal)
to $e$, the change in $\sN_G(X,\fe \mapsto e)$ will be
\begin{equation}
\label{eq:heuristic}
    \sN_G(X, \fe \mapsto e, e') \approx X(e'') \cdot | \cE_\triangle(G, \fe)| \cdot 2 p^{\se-3} n^{\sv-3}
    + |\cE_\wedge(G, \fe)| \cdot 2 p^{\se-2} n^{\sv-3}.
\end{equation}
On the other hand, the change in $\sN_\triangle(X,e)$ upon flipping $e'$ will be $6 X(e'')$
and the change in $\sN_\wedge(X,e)$ will be $2$, accounting for the $6$ automorphisms of the triangle and the $2$
automorphisms of the wedge.
Thus, we can match the changes of $\sN_G(X_t, \fe \mapsto e)$ along the Glauber dynamics with those of
\begin{equation}
\label{eq:matcher}
    \sC^\triangle_G(\fe) \cdot p^{\se-3} n^{\sv-3} \cdot \sN_\triangle(X_t,e)
    + \sC^\wedge_G(\fe) \cdot p^{\se-2} n^{\sv-3} \cdot \sN_\wedge(X_t,e),
\end{equation}
where $\sC^\triangle_G(\fe) = \frac{1}{3} |\cE_\triangle(G,\fe)|$ and $\sC^\wedge_G(\fe) = |\cE_\wedge(G,\fe)|$.
We remark that this definition shows that we have $C^\triangle_G(\fe), C^\wedge_G(\fe) \geq 0$ and that at least one of them is strictly
positive if $\fe$ has any neighboring edges in $G$, as required.
In any case, matching the changes of $\sN_G(X_t, \fe \mapsto e)$ along the Glauber dynamics with those of \eqref{eq:matcher}
will allow us to apply Theorem \ref{thm:barbour} to $\hsN_G(X, \fe \mapsto e)$ as defined in \eqref{eq:hsngxeee_def}
and obtain \eqref{eq:hajek_fe}.

Let us now prove the two aforementioned lemmas, after which we will make the above heuristic rigorous and prove
Proposition \ref{prop:hajek}.

\begin{proof}[Proof of Lemmas \ref{lem:stepbound_triangle} and \ref{lem:stepbound_wedge}]
Let us define
\begin{equation}
\label{eq:fchoice}
    f_\triangle(x) = \sN_G(x^{+e''}, \fe \mapsto e, \fe' \mapsto e')
    \qquad \text{or} \qquad
    f_\wedge(x) = \sN_G(x, \fe \mapsto e, \fe' \mapsto e').
\end{equation}
We will analyze $f_\triangle(x)$ under the hypotheses of Lemma \ref{lem:stepbound_triangle}
(i.e.\ if $\fe' \in \cE_\triangle(G,\fe)$) and analyze $f_\wedge(x)$ under the hypotheses
of Lemma \ref{lem:stepbound_wedge} (i.e.\ if $\fe' \in \cE_\wedge(G,\fe)$).
Shortly we will show that for any $\alpha \in (0,\frac{1}{2})$ there is some $c > 0$ such that
\begin{equation}
\label{eq:stepbound_concentration}
    \P[|f_*(X) - \E[f_*(X)]| > \gamma n^{\sv-3.5+\alpha}] \leq \Exp{-c n^{2\alpha}}
\end{equation}
for both $* = \triangle$ and $* = \wedge$ under the relevant hypotheses, which will finish the proof after we approximate
$\E[f_\triangle(X)]$ and $\E[f_\wedge(X)]$ closely enough, which we do now before proving \eqref{eq:stepbound_concentration}.
The approximation of these expectations follows a similar logic to the proof of Lemma \ref{lem:multilinear_applied},
and the proof of \eqref{eq:stepbound_concentration} is similar to the first part of the proof of Lemma \ref{lem:concentration_moments}.

First, if $\fe' \in \cE_\triangle(G,\fe)$, then $f_\triangle(x)$ may be written as a sum of indicators that various maps
are homomorphisms, one for each map $\cV(G) \to [n]$ which takes $\fe$ to $e$ and $\fe'$ to $e'$.
The number of such maps is $2 n^{\sv-3}$, the factor of $2$ coming from the choice of orientation of the
map $\fe \mapsto e$ (which determines the orientation of the map $\fe' \mapsto e'$), and the factor of $n^{\sv-3}$
simply accounting for the choice of where to map all remaining $\sv-3$ vertices of $G$.
All but $O(n^{\sv-4})$ of these maps are injective, which means that the indicator that the map is a homomorphism in $x$
can be written as a product of $\se-3$ distinct edge variables, recalling that $e''$ is included automatically in the definition
of $f_\triangle(x)$.
Thus we may apply Proposition \ref{prop:multilinear} and we find that for $X \sim \cM(n,p)$,
\begin{equation}
\label{eq:exptri}
    \E[f_\triangle(X)] = 2 n^{\sv-3} p^{\se-3} + O(n^{\sv-4}).
\end{equation}
Similar reasoning implies that if $\fe' \in \cE_\wedge(G,\fe)$ then
\begin{equation}
\label{eq:expwed}
    \E[f_\wedge(X)] = 2 n^{\sv-3} p^{\se-2} + O(n^{\sv-4}),
\end{equation}
the difference being that only $2$ edges are included for free instead of $3$.

Now let us prove \eqref{eq:stepbound_concentration}.
For this we will apply Theorem \ref{thm:lipschitzconcentration}, meaning we must calculate a bound on a Lipschitz vector $\cL$
for both $f_\triangle$ and $f_\wedge$.
In other words, we must bound the changes
\begin{equation}
    \partial_{e'''} f_*(x) = f_*(x^{+e'''}) - f_*(x^{-e'''})
\end{equation}
for any $e''' \in \edgeset$ and any $x \in \{0,1\}^{\edgeset}$, and for both $* = \triangle$ and $* = \wedge$.

First, if $e'''$ is either $e$ or $e'$ then $\partial_{e'''} f_*(x) = 0$ since neither function depends on the state of 
these two edges as they are always included in $x$ when doing the calculation.
Similarly, if $e'''$ is the third edge $e''$ of the triangle formed by $e$ and $e'$, then $\partial_{e'''} f_\triangle(x) = 0$.
On the other hand, $\partial_{e''} f_\wedge(x)$ may not be zero, but if $\fe' \in \cE_\wedge(G,\fe)$
then this always $\lesssim n^{\sv-4}$ since any map which covers $e''$ must be non-injective, since $\fe$ and $\fe'$
do not span a triangle in $G$.

Next, if $e'''$ is adjacent to the triangle spanned by $e$ and $e'$ but not a part of it, then
$\partial_{e'''} f_*(x) \lesssim n^{\sv-4}$ in both cases since there are $4$ vertices already accounted for by $e, e'$, and $e'''$
in any configuration (either these three edges form a $3$-star or a length-$3$ path).
Finally, if $e'''$ is not adjacent to the triangle spanned by $e$ and $e'$, then $\partial_{e'''} f_*(x) \lesssim n^{\sv-5}$
since there are three vertices accounted for in the triangle and two more in $e'''$.

Now, since there are $O(n)$ edges adjacent to the triangle spanned by $e$ and $e'$, in both cases the Lipschitz vector satisfies
\begin{equation}
    \| \cL \|_\infty \lesssim n^{\sv-4} \qquad \text{and} \qquad
    \| \cL \|_1 \lesssim n \cdot n^{\sv-4} + n^2 \cdot n^{\sv-5} \lesssim n^{\sv-3}.
\end{equation}
Therefore, by Theorem \ref{thm:lipschitzconcentration}, there are some constants $C, c > 0$
such that for all $\lambda \leq c n^{\sv-3}$ we have
\begin{equation}
    \P \left[
        \left|
            f_*(X) - \E[f_*(X)]
        \right| > \lambda
    \right] \leq 2 \Exp{- \frac{c\lambda^2}{n^{2\sv-7}} + e^{C\lambda n^{4-\sv} - c n}} + e^{-c n}
\end{equation}
Thus, plugging in $\lambda = \gamma n^{\sv-3.5+\alpha}$ we obtain
\begin{equation}
    \P \left[
        \left|
            f_*(X) - \E[f_*(X)]
        \right| > \gamma n^{\sv-3.5+\alpha}
    \right] \leq \Exp{- c \gamma^2 n^{2\alpha} + e^{C \gamma n^{\alpha + 0.5} - c n}} + e^{-cn}.
\end{equation}
Since $\alpha < \frac{1}{2}$, this finishes the proof of \eqref{eq:stepbound_concentration},
and by adjusting the constant and applying \eqref{eq:exptri} and \eqref{eq:expwed} we finish the proof
of Lemmas \ref{lem:stepbound_triangle} and \ref{lem:stepbound_wedge}, since
$n^{\sv-4} \ll n^{\sv-3.5+\alpha}$.
\end{proof}

With Lemmas \ref{lem:stepbound_triangle} and \ref{lem:stepbound_wedge} in hand, we turn to the proof of Proposition \ref{prop:hajek}.
As discussed above, it suffices to show \eqref{eq:hajek_fe}, namely that for any $\zeta \in (0, \frac{3}{4})$ and any $\gamma > 0$
there is some $c > 0$ such that for any $\fe \in \cE(G)$ we have
\begin{equation}
\label{eq:hajek_fe_again}
    \P\left[
        \left|
            \hsN_G(X,\fe \mapsto e) - \E[\hsN_G(X,\fe \mapsto e)]
        \right| > \gamma n^{\sv - 2.75 + \zeta}
    \right] < \Exp{- c n^{\frac{4}{3} \zeta}},
\end{equation}
where we recall from \eqref{eq:hsngxeee_def} that
\begin{equation}
\label{eq:hsngxeee_again}
    \hsN_G(x,\fe \mapsto e) = \sN_G(x,\fe \mapsto e)
    - \sC_G^\triangle(\fe) \cdot p^{\se-3} n^{\sv-3} \cdot \sN_\triangle(x,e)
    - \sC_G^\wedge(\fe) \cdot p^{\se-2} n^{\sv-3} \cdot \sN_\wedge(x,e)
\end{equation}
where, in light of \eqref{eq:heuristic},
$\sC_G^\triangle(\fe) = \frac{1}{3} |\cE_\triangle(G,\fe)|$ and $\sC_G^\wedge(\fe) = | \cE_\wedge(G,\fe) |$,
recalling that $\cE_\triangle(G,\fe)$ and $\cE_\wedge(G,\fe)$ are the sets of edges $\fe' \in \cE(G)$ which
share a vertex with $\fe$ and either span a triangle or \emph{do not} span a triangle with $\fe$ in $G$, respectively.

\begin{proof}[Proof of Proposition \ref{prop:hajek}]
As discussed above, we will fix $\fe \in \cE(G)$ and prove \eqref{eq:hajek_fe_again}.
For some $\alpha \in (0,\frac{1}{2})$ to be specified later,
let us define two sets of states where the concentration in the conclusion of Lemmas \ref{lem:stepbound_triangle}
and \ref{lem:stepbound_wedge} hold with this choice of $\alpha$ and with $\gamma = 1$:
\begin{align}
    \Pi_\triangle^\alpha &\coloneqq \left\{
        x \in \pball : 
        \begin{array}{c}
            \left| \sN_G(x^{+e''}, \fe \mapsto e, \fe' \mapsto e') - 2 p^{\se-3} n^{\sv-3} \right| \leq n^{\sv-3.5+\alpha}
            \quad \text{ for all } \fe' \in \cE_\triangle(G,\fe) \text{ and } \\
            \text{ all } e' \in \edgeset \setminus \{e\} \text{ with } e \cap e' \neq \emptyset \text{ and }
            e'' \text{ completing a triangle with } e,e' \text{ in } K_n
        \end{array}
    \right\}, \\
    \Pi_\wedge^\alpha &\coloneqq \left\{
        x \in \pball :
        \begin{array}{c}
            \left| \sN_G(x, \fe \mapsto e, \fe' \mapsto e') - 2 p^{\se-2} n^{\sv-3} \right| \leq n^{\sv-3.5+\alpha}
            \quad \text{ for}\\
            \text{all } \fe' \in \cE_\wedge(G,\fe)
            \text{ and all } e' \in \edgeset \setminus \{e\} \text{ with } e \cap e' \neq \emptyset
        \end{array}
    \right\},
\end{align}
and set $\Pi^\alpha = \Pi_\triangle^\alpha \cap \Pi_\wedge^\alpha$.
By Lemmas \ref{lem:stepbound_triangle} and \ref{lem:stepbound_wedge} and a union bound over the polynomially many edges $\fe'$ and $e'$,
we find that $\P[X \in \Pi^\alpha] \geq 1 - \Exp{-\Omega(n^{2\alpha})}$.
So let us apply Lemma \ref{lem:goodset_from_bigset} which results in a set $\Lambda^\alpha \sse \Lambda^*$ such that
\begin{equation}
\label{eq:stayinpi}
    \P[X_t^x \in \Pi^\alpha \text{ for all } 0 \leq t \leq n^3] \geq 1 - \Exp{-\Omega(n^{2\alpha})}
\end{equation}
for all $x$ with $\dh(x, \Lambda^\alpha) \leq 1$, and there is also some $z \in \Lambda^\alpha$ for which
\begin{equation}
\label{eq:stayinlambda}
    \P[X_t^z \in \Lambda^\alpha \text{ for all } 0 \leq t \leq n^3] \geq 1 - \Exp{-\Omega(n^{2\alpha})}.
\end{equation}
We will apply Theorem \ref{thm:barbour} with $\Lambda = \Lambda^\alpha$ and $f(x) = \hsN_G(x,\fe \mapsto e)$, with the choice $T = n^3$.
This means we need to validate the five conditions \ref{cond:variance}-\ref{cond:fbound} with appropriate constants
$V,J,\eps,\delta$, and $M$, where $V$ bounds a proxy for the variance of $f$, $J$ bounds a proxy for jumps of $f$ along the chain,
$\eps$ bounds the probability that the chain escapes $\Lambda^\alpha$, $\delta$ bounds the total variation distance of the chain to stationarity
at time $n^3$,
and $M$ is an overall bound on the difference between $f$ and its mean.

By the above construction, we may take $\epsilon = \Exp{- \Omega(n^{2\alpha})}$ since this is the probability the chain leaves $\Lambda^\alpha$
when started from $z$. 
Additionally, since $z \in \Lambda^*$, we may take $\delta = \Exp{-\Omega(n)}$ by the total variation distance bound in Lemma \ref{lem:mixing}.
And, since $|\hsN_G(x,\fe \mapsto e)| \lesssim n^{\sv-2}$ uniformly for all $x$, we may take $M = O(n^{\sv-2})$.

It just remains to bound $V$ and $J$ in conditions \ref{cond:variance} and \ref{cond:smalljumps}.
In preparation, let us consider a graph $x \in \Pi^\alpha$ and a single edge $e' \in \edgeset \setminus \{e\}$.
If $e'$ does not share a vertex with $e$, then by \eqref{eq:notneighbors} and the fact that $\sN_\triangle(x,e)$ and $\sN_\wedge(x,e)$
do not change when the status of $e'$ is flipped, we have
\begin{equation}
    \hsN_G(x^{+e'}, \fe \mapsto e) - \hsN_G(x^{-e'}, \fe \mapsto e) = \sN_G(x, \fe \mapsto e, e') \lesssim n^{\sv-4}.
\end{equation}
On the other hand, if $e'$ does share a vertex with $e$, and $e''$ completes the triangle in $K_n$, then by
\eqref{eq:Ngxeeerepresentation} and the subsequent discussion,
recalling our definitions of $\sC^\triangle_G(\fe) = \frac{1}{3} |\cE_\triangle(G,\fe)|$ and
$\sC^\wedge_G(\fe) = |\cE_\wedge(G,\fe)|$, we have
\begin{align}
    &\hsN_G(x^{+e'}, \fe \mapsto e) - \hsN_G(x^{-e'}, \fe \mapsto e) \\
    &\qquad = \sum_{\fe' \in \cE_\triangle(G,\fe)} \left(
        x(e'') \cdot \sN_G(x^{+e''}, \fe \mapsto e, \fe' \mapsto e')
        - \frac{1}{3} \cdot p^{\se - 3} n^{\sv - 3} \cdot \partial_{e'} \sN_\triangle(x,e)
    \right) \\
    &\qquad \qquad + \sum_{\fe' \in \cE_\wedge(G,\fe)} \left(
        \sN_G(x, \fe \mapsto e, \fe' \mapsto e')
        - p^{\se - 2} n^{\sv - 3} \cdot \partial_{e'} \sN_\wedge(x,e)
    \right) \\
    &\qquad \qquad + O(n^{\sv-4}).
\end{align}
Now recalling that a triangle has $6$ automorphisms and a wedge has $2$, we have $\partial_{e'} \sN_\triangle(x,e) = 6 x(e'')$
and $\partial_{e'} \sN_\wedge(x,e) = 2$.
Therefore the above expression is equal to
\begin{align}
    &x(e'') \cdot \sum_{\fe' \in \cE_\triangle(G,\fe)} \left(
        \sN_G(x^{+e''}, \fe \mapsto e, \fe' \mapsto e')
        - 2 p^{\se - 3} n^{\sv - 3}
    \right) \\
    &\qquad + \sum_{\fe' \in \cE_\wedge(G,\fe)} \left(
        \sN_G(x, \fe \mapsto e, \fe' \mapsto e')
        - 2 p^{\se - 2} n^{\sv - 3}
    \right) \\
    &\qquad + O(n^{\sv-4}).
\end{align}
Thus, since $x \in \Pi^\alpha$, the above expression is $\lesssim n^{\sv-3.5+\alpha}$, since $\Pi^\alpha$ is exactly where
all of the (finitely many) terms in the sums above are $\lesssim n^{\sv-3.5+\alpha}$.
Summarizing, if $x \in \Pi^\alpha$ then we have the following bounds:
\begin{equation}
\label{eq:stepbound}
    \hsN_G(x^{+e'}, \fe \mapsto e) - \hsN_G(x^{-e'}, \fe \mapsto e)
    \lesssim \begin{cases}
        n^{\sv-4} & \text{if } e \cap e' = \emptyset, \\
        n^{\sv-3.5+\alpha} & \text{if } e \cap e' \neq \emptyset.
    \end{cases}
\end{equation}
Note also that the expression on the left-hand side above is nonnegative.

Now notice that $\Pi^\alpha$ is an interval in the sense that if $a, b \in \Pi^\alpha$ and $c$ is such that with
$a \preceq c \preceq b$, then $c \in \Pi^\alpha$ as well.
Here $\preceq$ denotes the partial order of the Hamming cube where $a \preceq b$ if the edges in $a$ form a subset of the
edges in $b$.
This interval property holds because $\Pi^\alpha$ may be expressed as the intersection of many inverse images of intervals under increasing functions
$\{ 0, 1 \}^{\edgeset} \to \R$.
Therefore, if $a, b \in \Pi^\alpha$ with $a \preceq b$ then we may add one edge at a time, starting from $a$ and ending at $b$,
all while staying within $\Pi^\alpha$.
Applying \eqref{eq:stepbound} at each step in this process, we obtain
\begin{equation}
\label{eq:intervalbound}
    \left|
        \hsN_G(a, \fe \mapsto e) - \hsN_G(b, \fe \mapsto e)
    \right| \lesssim n^{\sv-3.5+\alpha} \cdot \dloce(a,b) + n^{\sv-4} \cdot \dh(a,b),
\end{equation}
where we recall from \eqref{eq:dlocdef} that
\begin{equation}
    \dloce(x,x') \coloneqq \sum_{e' \cap e \neq \emptyset} |x(e') - x'(e')|.
\end{equation}
We will use \eqref{eq:intervalbound} to provide the bounds on both of $V$ and $J$
in conditions \ref{cond:variance} and \ref{cond:smalljumps}.

Let us begin with condition \ref{cond:smalljumps}, for which we need to bound
\begin{equation}
    \left|
        \E \left[\hsN_G(X_t^x, \fe \mapsto e)\right]
        - \E \left[\hsN_G(X_t^{x'}, \fe \mapsto e)\right]
    \right|
\end{equation}
for all $t \geq 0$, whenever $x \in \Lambda^\alpha$ and $\dh(x,x') \leq 1$.
Since $\dloce(a,b) \leq \dh(a,b)$ and $n^{\sv-3.5+\alpha} \gg n^{\sv-4}$, we may simplify the bound \eqref{eq:intervalbound},
and obtain
\begin{equation}
\label{eq:simpleinterval}
    \left| \hsN_G(a, \fe \mapsto e) - \hsN_G(b, \fe \mapsto e) \right| \lesssim n^{\sv-3.5+\alpha} \dh(a,b)
\end{equation}
as long as both $a,b$ are in $\Pi^\alpha$ and $a \preceq b$ or vice versa.
Now, for $t \leq n^3$ the probability that either $X_t^x$ or $X_t^{x'}$ escapes $\Pi^\alpha$
is $\Exp{-\Omega(n^{2\alpha})}$ by \eqref{eq:stayinpi}.
Additionally, since $\Lambda^\alpha \subseteq \Lambda^*$, by Proposition \ref{prop:goodset} both chains have probability at least
$1 - e^{-\Omega(n)}$ to stay within $\phalfball$ up to time $n^3$; if this is the case, then the ordering inherited from $x \preceq x'$
or $x' \preceq x$ (which holds because $\dh(x,x') \leq 1$)
will be preserved, i.e.\ we will have $X_t^x \preceq X_t^{x'}$ or vice versa.
This is because we are using the monotone coupling which preserves monotonicity, except possibly for when the chains
attempt to leave the metastable well $\pball$, which cannot happen if they remain in $\phalfball$.

In any case, for $t \leq n^3$, by the reasoning in the above paragraph,
the second bound in Lemma \ref{lem:mixing} yields the existence of some constant $\kappa > 0$ for which we have
\begin{align}
\label{eq:linfbnd}
    \E \left[\left| \hsN_G(X_t^x, \fe \mapsto e) - \hsN_G(X_t^{x'}, \fe \mapsto e) \right| \right]
    &\lesssim n^{\sv-3.5+\alpha} \left( 1 - \frac{\kappa}{n^2} \right)^t + \Exp{-\Omega(n^{2\alpha})} \\
    &\lesssim n^{\sv-3.5+\alpha}.
\end{align}
On the other hand, again by the second bound in Lemma \ref{lem:mixing} and Markov's inequality,
the chains $X_t^x$ and $X_t^{x'}$ have probability at least $1 - e^{-\Omega(n)}$ of meeting before time $n^3$.
So, since
\begin{equation}
    \left| \hsN_G(X_t^x, \fe \mapsto e) - \hsN_G(X_t^{x'}, \fe \mapsto e) \right|
\end{equation}
is deterministically polynomially bounded and is zero with probability at least $1 - \Exp{-\Omega(n)}$ for $t > n^3$,
we have an exponentially small bound on the expectation of this quantity whenever $t > n^3$.
All of this means that we can take $J \lesssim n^{\sv - 3.5 + \alpha}$ in condition \ref{cond:smalljumps}.

Finally we turn to condition \ref{cond:variance}.
We need to come up with a suitable value of $V$ which bounds the following variance proxy for all $x \in \Lambda^\alpha$:
\begin{align}
    &\sum_{x' \in \pball} P(x,x') \sum_{t=0}^{n^3-1} \left(
        \E \left[\hsN_G(X_t^x, \fe \mapsto e)\right]
        - \E \left[\hsN_G(X_t^{x'}, \fe \mapsto e)\right]
    \right)^2 \\
    &\qquad \leq \frac{1}{\binom{n}{2}} \sum_{e' \in \edgeset}
    \sum_{t=0}^{n^3-1} \left(
        \E \left[\hsN_G(X_t^x, \fe \mapsto e)\right]
        - \E \left[\hsN_G(X_t^{x^{\oplus e'}}, \fe \mapsto e)\right]
    \right)^2,
\label{eq:goalsum}
\end{align}
where we use the notation $x^{\oplus e'}$ as introduced in Section \ref{sec:inputs_consequences_concentration}
to denote the graph which is the same as $x$ but with the status of $e'$ flipped.

Let us first consider the terms of \eqref{eq:goalsum} corresponding to $e'$ with $e \cap e' \neq \emptyset$.
For these, we may apply the simplified bound \eqref{eq:simpleinterval}
and then the second bound in Lemma \ref{lem:mixing}, as in the derivation of \eqref{eq:linfbnd} above, to obtain
the following bound on the inner sum:
\begin{align}
    &\sum_{t=0}^{n^3-1} \left(
        \E \left[ \hsN_G(X_t^x, \fe \mapsto e) \right]
        - \E \left[ \hsN_G(X_t^{x^{\oplus e'}}, \fe \mapsto e) \right]
    \right)^2 \\
    &\qquad \lesssim n^{2\sv - 7 + 2\alpha} \sum_{t=0}^{n^3-1} \left( \E \left[ \dh(X_t^x, X_t^{x^{\oplus e'}}) \right] \right)^2 + \Exp{-\Omega(n^{2\alpha})} \\
    &\qquad \lesssim n^{2\sv - 7 + 2\alpha} \sum_{t=0}^{n^3-1} \left( 1 - \frac{\kappa}{n^2} \right)^{2t} + \Exp{-\Omega(n^{2\alpha})} \\
    &\qquad \lesssim n^{2\sv - 5 + 2\alpha},
\end{align}
in the case where $e \cap e' \neq \emptyset$.

Now let us consider the terms of \eqref{eq:goalsum} corresponding to $e'$ with $e \cap e' = \emptyset$.
Going back to the unsimplified bound \eqref{eq:intervalbound}, and again using the fact that $X_t^x, X_t^{x^{\oplus e'}} \in \Pi^\alpha$
and are monotonically ordered with probability at least $1 - \Exp{- \Omega(n^{2\alpha})}$, we have
\begin{align}
    &\left|
        \E \left[ \hsN_G(X_t^x, \fe \mapsto e) \right]
        - \E \left[ \hsN_G(X_t^{x^{\oplus e'}}, \fe \mapsto e) \right]
    \right| \\
    &\qquad \lesssim n^{\sv-3.5+\alpha} \cdot \E \left[ \dloce(X_t^x, X_t^{x^{\oplus e'}})\right]
    + n^{\sv-4} \cdot \E \left[ \dh(X_t^x, X_t^{x^{\oplus e'}})\right]
    + \Exp{-\Omega(n^{2\alpha})}.
\end{align}
Now since $e'$ is not adjacent to $e$, by Corollary \ref{cor:dloce}, we have
\begin{equation}
    \E \left[ \dloce(X_t^x, X_t^{x^{\oplus e'}}) \right] \lesssim \frac{1}{n}.
\end{equation}
Now since $\E \left[ \dh(X_t^x, X_t^{x^{\oplus e'}}) \right] \lesssim 1$ by the second bound in Lemma \ref{lem:mixing},
and since $\alpha < \frac{1}{2}$, we have
\begin{equation}
\label{eq:l1bnd}
    \left|
        \E \left[ \hsN_G(X_t^x, \fe \mapsto e) \right]
        - \E \left[ \hsN_G(X_t^{x^{\oplus e'}}, \fe \mapsto e) \right]
    \right| \lesssim n^{\sv-4}.
\end{equation}
Thus by peeling off one factor of the square in the inner sum of \eqref{eq:goalsum} and applying \eqref{eq:l1bnd},
then applying \eqref{eq:linfbnd} to the other factor inside the sum,
and finally applying the second bound of Lemma \ref{lem:mixing} again,
we obtain
\begin{align}
    &\sum_{t=0}^{n^3-1} \left(
        \E \left[ \hsN_G(X_t^x, \fe \mapsto e) \right]
        - \E \left[ \hsN_G(X_t^{x^{\oplus e'}}, \fe \mapsto e) \right]
    \right)^2 \\
    &\qquad \lesssim n^{\sv-4} \sum_{t=0}^{n^3-1} \left|
        \E \left[\hsN_G(X_t^x, \fe \mapsto e)\right]
        -\E \left[\hsN_G(X_t^{x^{\oplus e'}}, \fe \mapsto e)\right]
    \right| \\
    &\qquad \lesssim n^{\sv-4} \sum_{t=0}^{n^3-1} n^{\sv-3.5+\alpha} \cdot \E \left[ \dh(X_t^x, X_t^{x^{\oplus e'}})\right] + \Exp{-\Omega(n^{2\alpha})} \\
    &\qquad \lesssim n^{2\sv-7.5+\alpha} \sum_{t=0}^{n^3-1} \left(1 + \frac{\kappa}{n^2} \right)^t + \Exp{-\Omega(n^{2\alpha})} \\
    &\qquad \lesssim n^{2\sv-5.5+\alpha}.
\end{align}
Finally, combining the analysis of the cases where $e' \cap e \neq \emptyset$ and where $e' \cap e = \emptyset$,
and also using the fact that $\alpha < \frac{1}{2}$, we find that the original sum \eqref{eq:goalsum} we needed
to bound is
\begin{equation}
    \lesssim \frac{1}{n^2} \left( n \cdot n^{2\sv - 5 + 2\alpha} + n^2 \cdot n^{2\sv-5.5+\alpha} \right)
    = n^{2\sv-6+2\alpha} + n^{2\sv-5.5+\alpha}
    \lesssim n^{2\sv-5.5+\alpha},
\end{equation}
using the fact that $\alpha < \frac{1}{2}$ for the last inequality.
Therefore condition \ref{cond:variance} holds with $V \lesssim n^{2\sv-5.5+\alpha}$.

In summary, all conditions of Theorem \ref{thm:barbour} hold with the following choices:
\begin{equation}
    V = C n^{2\sv-5.5+\alpha}, \quad
    J = C n^{\sv-3.5+\alpha}, \quad
    \epsilon = \Exp{-c n^{2\alpha}}, \quad
    \delta = \Exp{-c n}, \quad
    M = C n^{\sv-2}
\end{equation}
for some constants $C$ and $c$.
Thus the conclusion of Theorem \ref{thm:barbour} yields
\begin{align}
    &\P \left[ \left| \hsN_G(X, \fe \mapsto e) - \E \left[ \hsN_G(X, \fe \mapsto e) \right] \right| > \lambda \right] \\
    &\qquad \leq 2 \Exp{- \frac{(\lambda - 2 C n^{\sv-2} e^{-c n})^2}{2 C n^{2\sv-5.5+\alpha} + \frac{4}{3} C n^{\sv-3.5+\alpha} (\lambda - 2 C n^{\sv-2} e^{-cn})}}
    + e^{- c n^{2\alpha}} + e^{-cn}.
\end{align}
Plugging in $\lambda = \gamma n^{\sv - 2.75 + \zeta}$, the second term in the denominator is dominated by the first exactly when
\begin{equation}
    2 \sv - 5.5 + \alpha \geq 2 \sv - 6.25 + \alpha + \zeta,
\end{equation}
which means we should take $\zeta \leq 0.75$ for this to hold.
Under this assumption, by simplifying the expression, changing the constants, and using the fact that $\alpha < \frac{1}{2}$, we find that
\begin{equation}
    \P \left[ \left| \hsN_G(X, \fe \mapsto e) - \E \left[ \hsN_G(X, \fe \mapsto e) \right] \right| > \gamma n^{\sv-2.75+\zeta} \right]
    \leq \Exp{- c n^{2\zeta-\alpha}} + \Exp{-c n^{2\alpha}}.
\end{equation}
The optimal choice of $\alpha$ here is thus $\alpha = \frac{2}{3} \zeta$, which is in $(0, \frac{1}{2})$ since $\zeta \in (0, \frac{3}{4})$.
This gives an upper bound of $\Exp{- c n^{\frac{4}{3} \zeta}}$ and finishes the proof
of \eqref{eq:hajek_fe_again}, which completes the proof of Proposition \ref{prop:hajek}.
\end{proof}

Now we use Proposition \ref{prop:hajek} and the decomposition of the Hamiltonian differential \eqref{eq:peh_expansion_3}
to provide an analogous approximation for the Hamiltonian differential itself.

\begin{corollary}
\label{cor:hajek}
There are constants $\sC^\triangle, \sC^\wedge \geq 0$,
with at least one of them strictly positive under the nondegeneracy assumption \eqref{eq:nondegen},
such that, if we define
\begin{equation}
    \heH(x) \coloneqq n^2 \partial_e \sH(x) - \frac{1}{n} \left( \sC^\triangle \cdot \sN_\triangle(x,e) + \sC^\wedge \cdot \sN_\wedge(x,e) \right),
\end{equation}
then for any $\zeta \in (0, \frac{3}{4})$ and any $\gamma > 0$ there is some $c > 0$ such that
\begin{equation}
    \P \left[
        \left|
            \heH(X) - \E \left[ \heH(X) \right]
        \right| > \gamma n^{-0.75+\zeta}
    \right] \leq \Exp{- c n^{\frac{4}{3} \zeta}}.
\end{equation}
\end{corollary}

\begin{proof}[Proof of Corollary \ref{cor:hajek}]
Let us define
\begin{equation}
    \sC^\triangle = \sum_{j=0}^K \beta_j p^{\se_j-3} \sC^\triangle_{G_j}
    \qquad \text{and} \qquad
    \sC^\wedge = \sum_{j=0}^K \beta_j p^{\se_j-2} \sC^\wedge_{G_j},
\end{equation}
where $\sC^\triangle_G$ and $\sC^\wedge_G$ are the constants guaranteed by Proposition \ref{prop:hajek} for any graph $G$,
and $G_j$ are the graphs in the ERGM specification.
Note that under the nondegeneracy assumption \eqref{eq:nondegen}, at least one of these graphs
is not a disjoint union of edges, and note also that $\sC^\triangle_{G_0} = \sC^\wedge_{G_0} = 0$ since $G_0$ is a single edge.
So, since $\beta_j > 0$ for $j \geq 1$, at least one of the above constants $\sC^\triangle, \sC^\wedge$ is strictly positive.
Now, recalling \eqref{eq:peh_expansion_3}, with these choices of $\sC^\triangle$ and $\sC^\wedge$, if
\begin{equation}
    \left| \heH(X) - \E[\heH(X)] \right| > \gamma n^{- 0.75 + \zeta},
\end{equation}
then for some $j > 0$ we must also have
\begin{equation}
\label{eq:nbig}
    \frac{\beta_j}{n^{\sv_j-2}} \left| \hsN_{G_j} (X) - \E[\hsN_{G_j}(X)] \right| > \frac{\gamma}{K} n^{- 0.75 + \zeta},
\end{equation}
noting that the $j=0$ term in the sum we are decomposing above is zero, since $\sN_{G_0}(x,e) \equiv 1$ for $G_0$ a single edge.
In any case, by Proposition \ref{prop:hajek}, the event \eqref{eq:nbig} happens with probability at most $\Exp{- c_j n^{\frac{4}{3} \zeta}}$
for some constant $c_j$ which may be different for different values of $j$.
Nevertheless, since there are only finitely many options for $j$, we immediately obtain the result.
\end{proof}

%% file: sections/hamiltonian/fkg.tex
\subsection{Fluctuation lower bounds via approximate FKG inequality}
\label{sec:hamiltonian_fkg}

In this section we use Corollary \ref{cor:hajek} as well as the approximate FKG inequality
of Theorem \ref{thm:fkg} to derive a lower bound on the fluctuations of $n^2 \partial_e \sH(X)$
for $X \sim \cM(n,p)$, which will be useful for proving the lower bound in our
main result, Theorem \ref{thm:wasserstein}.
We begin with a lower bound on the fluctuations of the quantity which Corollary \ref{cor:hajek}
compares $n^2 \partial_e \sH(X)$ with.

\begin{proposition}
\label{prop:nfluct}
Let $X \sim \cM(n,p)$ and define
For some $\sC^\triangle, \sC^\wedge \geq 0$, define
\begin{equation}
    \sS(x) = \sC^\triangle \cdot \sN_\triangle(x,e) + \sC^\wedge \cdot \sN_\wedge(x,e).
\end{equation}
If either $\sC^\triangle > 0$ or $\sC^\wedge > 0$ then there are some $\eps, \delta > 0$ such that, for $X \sim \cM(n,p)$, we have
\begin{equation}
    \P \left[
        |\sS(X) - \E[\sS(X)]| > \delta \sqrt{n}
    \right] > \eps.
\end{equation}
\end{proposition}

\begin{proof}[Proof of Proposition \ref{prop:nfluct}]
Since $\sC^\triangle, \sC^\wedge \geq 0$, we have
\begin{align}
    \Var[\sS(X)] &= (\sC^\triangle)^2 \Var[\sN_\triangle(X,e)]
        + (\sC^\wedge)^2 \Var[\sN_\wedge(X,e)]
        + \sC^\triangle \sC^\wedge \Cov[\sN_\triangle(X,e), \sN_\wedge(X,e)] \\
    &\geq (\sC^\triangle)^2 \Var[\sN_\triangle(X,e)]
        + (\sC^\wedge)^2 \Var[\sN_\wedge(X,e)] - e^{-\Omega(n)},
\end{align}
using Theorem \ref{thm:fkg} for the lower bound on the covariance.
Now, writing $e = \{u,v\}$, we also have
\begin{align}
    \Var[\sN_\triangle(X,e)] &= \sum_{w \notin \{ u,v \}} \Var[X(\{u,w\}) X(\{v,w\})] \\
    &\qquad + \sum_{w_1,w_2 \notin \{ u,v \}} \Cov[X(\{u,w_1\}) X(\{v,w_1\}), X(\{u,w_2\}) X(\{v,w_2\})] \\
    &\gtrsim n - n^2 e^{-\Omega(n)} \gtrsim n,
\end{align}
using Theorem \ref{thm:fkg} again as well as the fact that the variances in the first sum above
are all the same by symmetry.
Similarly, we also have $\Var[\sN_\wedge(X,e)] \gtrsim n$.
Thus, since at least one of $\sC^\triangle$ and $\sC^\wedge$ is strictly positive,
we find that $\Var[\sS(X)] \gtrsim n$ as well.

Next we will apply the concentration inequality of Theorem \ref{thm:lipschitzconcentration}
to \emph{upper} bound the fluctuations of $\sS(X)$ by $O(\sqrt{n})$, which will then allow us
to use the variance lower bound of $\Omega(n)$ for $\sS(X)$ to obtain a matching lower bound of $\Omega(\sqrt{n})$ on
the fluctuations themselves.

To apply Theorem \ref{thm:lipschitzconcentration} we will need to calculate a bound for the Lipschitz vector $\cL$ of $\sS(x)$.
First, since $\sS(x)$ does not depend on $x(e)$, we have $\cL_e = 0$.
Next, if $e'$ is adjacent to $e$ then changing $x(e')$ may change $\sN_\triangle(x,e)$
and $\sN_\wedge(x,e)$ by at most $O(1)$ each.
Finally, if $e'$ is not adjacent to $e$ then changing $x(e')$ does not change $\sS(x)$.
So we have
\begin{equation}
    \| \cL \|_\infty \lesssim 1
    \qquad \text{and} \qquad
    \| \cL \|_1 \lesssim n.
\end{equation}
Thus Theorem \ref{thm:lipschitzconcentration} yields $C,c > 0$ such that for all $\lambda \leq c n$
we have
\begin{equation}
\label{eq:sconc}
    \P \left[ |\sS(X) - \E[\sS(X)]| > \lambda \right]
    \leq 2 \Exp{- c \frac{\lambda^2}{n} + c e^{C \lambda - c n}} + e^{-c n}.
\end{equation}
We will use this result shortly.

First, let us set $T = (\sS(X) - \E[\sS(X)])^2$ and note that for any $M$ we have
\begin{equation}
\label{eq:varcompar}
    n \lesssim \Var[\sS(X)] = \E[T \ind{T \leq M n}] + \E[T \ind{T > M n}].
\end{equation}
Now let us bound the right-most term for large $M$ using \eqref{eq:sconc}.
We have
\begin{align}
    \E[T \ind{T > M n}] &= \int_0^\infty \P[T \ind{T > M n} \geq t] \,dt \\
    &\leq \int_0^{M n} \P[T \geq M n] \,dt + \int_{M n}^{\frac{c^2}{C^2} n^2} \P[T \geq t] \,dt +
    \int_{\frac{c^2}{C^2} n^2}^{O(n^2)} \P[T \geq t] \,dt,
\label{eq:threeintegrals}
\end{align}
where the $O(n^2)$ term in the right-most integral is some deterministic upper bound for $T$,
and the $C,c$ in the second and third integrals are the same as in the conclusion of
Theorem \ref{thm:lipschitzconcentration} and \eqref{eq:sconc} above.
Since $\sqrt{M n} \leq c n$ for $n$ large enough, we may
use \eqref{eq:sconc} with $\lambda = \sqrt{M n}$ to bound the first integral in \eqref{eq:threeintegrals} as
\begin{equation}
    \int_0^{M n} \P[T \geq M n] \,dt \leq M n \cdot \left( 2 \Exp{- c M + c e^{C \sqrt{M n} - c n}} + e^{- cn} \right)
    \leq 3 M n e^{- c M}
\end{equation}
for $n$ large enough.
Similarly, we may use \eqref{eq:sconc} with $\lambda = t$ to bound the second integral in \eqref{eq:threeintegrals} as
\begin{align}
    \int_{M n}^{\frac{c^2}{C^2} n^2} \P[T \geq t] \,dt
    &\leq \int_{M n}^{\frac{c^2}{C^2} n^2}
    \left( 2 \Exp{- c \frac{t}{n} + c e^{C \sqrt{t} - c n}} + e^{-c n} \right) \,dt \\
    &\leq 2 e^c \int_{M n}^{\frac{c^2}{C^2} n^2} e^{- c t / n} \,dt + \frac{c^2}{C^2} n^2 e^{- cn} \\
    &\leq 2 e^c \frac{n}{c} e^{- c M} + e^{- \Omega(n)}.
\end{align}
Finally, the third integral in \eqref{eq:threeintegrals} can be bounded as
\begin{align}
    \int_{\frac{c^2}{C^2} n^2}^{O(n^2)} \P[T \geq t] \,dt 
    &\leq \int_{\frac{c^2}{C^2} n^2}^{O(n^2)} \P\left[T \geq \frac{c^2}{C^2} n^2 \right] \,dt \\
    &\leq O(n^2) \cdot \left( 2 \Exp{- \frac{c^3}{C^2} n + c} + e^{- c n} \right) \\
    &\leq e^{-\Omega(n)}.
\end{align}
Thus we find that
\begin{equation}
    \E[T \ind{T > M n}] \leq \left(3 M + \frac{2 e^c}{c}\right) e^{- cM} n + e^{-\Omega(n)}.
\end{equation}
By choosing $M$ large enough
we may assure that $\E[T \ind{T > M n}]$ is at most an arbitrarily small multiple of $n$, and so,
recalling \eqref{eq:varcompar}, we have
\begin{equation}
    \E[T \ind{T \leq M n}] \gtrsim n.
\end{equation}
Since $T \ind{T \leq M n} \leq M n$,
there must be some $\delta, \eps > 0$ such that
we have $T \geq \delta^2 n$ with probability
at least $\eps$.
This finishes the proof.
\end{proof}

With Proposition \ref{prop:nfluct} in hand, we can prove the main ingredient for the
lower bound on the Wasserstein distance in Theorem \ref{thm:wasserstein}, which is a useful form
of the lower bound on the fluctuations of the Hamiltonian differential.
This is a matching lower bound for the $k=1$ case of Corollary \ref{cor:concentration_marginal}.

\begin{corollary}
\label{cor:dhfluct}
For $X \sim \cM(n,p)$, under the nondegeneracy assumption \eqref{eq:nondegen} we have
\begin{equation}
    \E \left[ \left| n^2 \partial_e \sH(X) - 2 \sH'(p) \right| \right] \gtrsim n^{-1/2}.
\end{equation}
\end{corollary}

\begin{proof}[Proof of Corollary \ref{cor:dhfluct}]
First note that by Theorem \ref{thm:marginal} and Lemma \ref{lem:multilinear_applied}, we have
\begin{equation}
    2 \sH'(p) = 2 \sH'(\E[X(e)]) + O(n^{-1}) = \E[n^2 \partial_e \sH(X)] + O(n^{-1}).
\end{equation}
We also used the fact that $\sH'$ is Lipschitz in the above derivation.
So it suffices to show that
\begin{equation}
    \E \left[ \left| n^2 \partial_e \sH(X) - \E[n^2 \partial_e \sH(X)] \right| \right] \gtrsim n^{-1/2}.
\end{equation}
Now recall the notation
\begin{equation}
    \sS(x) = \sC^\triangle \cdot \sN_\triangle(x,e) + \sC^\wedge \cdot \sN_\wedge(x,e)
\end{equation}
from the statement of Proposition \ref{prop:nfluct}, where we now plug in the valuse of $\sC^\triangle$
and $\sC^\wedge$ given by Corollary \ref{cor:hajek}.
Corollary \ref{cor:hajek} states that for all $\zeta \in (0,\frac{3}{4})$ and $\gamma > 0$ there is some $c > 0$
such that
\begin{equation}
    \P \left[
        \left| \left( n^2 \partial_e \sH(X) - \E[n^2 \partial_e \sH(X)] \right)
        - \frac{1}{n} (\sS(X) - \E[\sS(X)])
        \right| > \gamma n^{-\frac{3}{4} + \zeta} \right] \leq \Exp{- c n^{\frac{4}{3} \zeta}}.
\end{equation}
Taking $\zeta = \frac{1}{8}$ and $\gamma = 1$ for instance, we see that
\begin{equation}
    n^2 \partial_e \sH(X) - \E[n^2 \partial_e \sH(X)]
    \qquad \text{and} \qquad
    \frac{1}{n} (\sS(X) - \E[\sS(X)])
\end{equation}
are within $n^{- \frac{5}{8}}$ of each other with probability at least $1 - e^{- c n^{\frac{1}{6}}}$.
Now by Proposition \ref{prop:nfluct} (since the nondegeneracy assumption \eqref{eq:nondegen}
ensures that $\sC^\triangle>0$ or $\sC^\wedge > 0$), we have
\begin{equation}
    \P \left[ \frac{1}{n} |\sS(X) - \E[\sS(X)]| > \delta n^{-1/2} \right] > \eps,
\end{equation}
and so we also have
\begin{equation}
    \P \left[ \left| n^2 \partial_e \sH(X) - \E[n^2 \partial_e \sH(X)]\right|
    > \delta n^{-\frac{1}{2}} - n^{-\frac{5}{8}}\right] > \eps - e^{- c n^{\frac{1}{6}}}.
\end{equation}
This finishes the proof.
\end{proof}

Before moving on, we remark that the reduction to $\sN_\triangle(X,e)$ and $\sN_\wedge(X,e)$ given by
Proposition \ref{prop:hajek} is necessary to obtain the bound given by Corollary \ref{cor:dhfluct}.
Indeed, attempting to apply a strategy like the one presented in Proposition \ref{prop:nfluct} directly to
$n^2 \partial_e \sH(X)$ itself will fail in general.
For instance, if $K=1$ and $G_1 = \square$ is a square, then letting $e = \{u,v\}$ we find that 
the variance of $\sN_\square(X,e)$ is
\begin{align}
    & \sum_{w_1,w_2} \Var \left[ X(\{u,w_1\}) X(\{v,w_2\}) X(\{w_1,w_2\}) \right] \\
    &\qquad + \sum_{w_1, w_2, w_3, w_4} \Cov \left[
        X(\{u,w_1\}) X(\{v,w_2\}) X(\{w_1,w_2\}),
        X(\{u,w_3\}) X(\{v,w_4\}) X(\{w_3,w_4\})
    \right].
\end{align}
Now the first sum is only $\gtrsim n^2$ while we know from Propositions \ref{prop:hajek}
and \ref{prop:nfluct} that the true variance is of order $n^3$, since $\sN_\square(X,e)$ is well-approximated
by a constant times $n \cdot \sN_\wedge(X,e)$ (we have $\sC^\triangle = 0$ in this particular case).
Thus the second sum above must contain the main contribution to the variance, but a priori we have no way
to \emph{lower bound} the covariances in that sum, other than by a negative number using Theorem \ref{thm:fkg}.

%% file: sections/wasserstein.tex
\section{Wasserstein distance}
\label{sec:wasserstein}

In this section we will prove our main result, Theorem \ref{thm:wasserstein}, giving matching upper and lower bounds
on the Wasserstein distance between $\cM(n,p)$ and $\cG(n,p)$.
First, in Section \ref{sec:wasserstein_coupling} we describe the optimal coupling between $\cM(n,p)$
and $\cG(n,p)$, and then in Section \ref{sec:wasserstein_proof} we use this description as well
as the preparation of the previous sections to give the proof.

\input{sections/wasserstein/coupling}
\input{sections/wasserstein/proof}
\input{sections/wasserstein/alternative}

%% file: sections/wasserstein/coupling.tex
\subsection{The optimal coupling}
\label{sec:wasserstein_coupling}

Recall from \eqref{eq:wasdef} that the Wasserstein distance $\dW(\cM(n,p),\cG(n,p))$ is defined as
\begin{equation}
\label{eq:wasdef_again}
    \E \left[ \dh(X,Y) \right] = \E \left[ \sum_{e \in \edgeset} \ind{X(e) \neq Y(e)} \right]
\end{equation}
under the coupling $(X,Y)$ with marginals $X \sim \cM(n,p)$ and $Y \sim \cG(n,p)$ which minimizes
the above expectation.
While there may be multiple optimal couplings which achieve this, they must all satisfy a
local optimality condition under resampling edges, described in the following lemma.
We state the lemma in a slightly more general context as it does not depend on the graph structure
or the particular details of the two measures $\cM(n,p)$ or $\cG(n,p)$.

\begin{lemma}
\label{lem:optimality}
Fix any Wasserstein-optimal coupling $(X,Y)$ between two measures $\mu$ and $\nu$
on $\{ 0, 1 \}^\cE$ for some finite set $\cE$.
For any fixed coordinate $e \in \cE$, the conditional distribution of $(X(e), Y(e))$
given $(X(e'),Y(e'))$ for all $e' \neq e$ must be that of a pair of optimally
coupled Bernoulli random variables with means
\begin{equation}
\label{eq:optimalprobs}
    \E_\mu[X(e) | X(e') \text{ for } e' \neq e ]
    \qquad \text{and} \qquad
    \E_\nu[Y(e) | Y(e') \text{ for } e' \neq e]
\end{equation}
respectively.
\end{lemma}

\begin{proof}[Proof of Lemma \ref{lem:optimality}]
Suppose this were not the case for some fixed $e \in \cE$, and let $(X,Y)$ denote a sample from
the supposedly optimal coupling.
Then define $(X',Y')$ by first sampling $(X,Y)$ and then resampling the $e$th coordinates of both
$X$ and $Y$ according to an optimal coupling of Bernoulli random variables with the probabilities
described in \eqref{eq:optimalprobs}.
Then since no other edges are resampled we have
\begin{align}
    \E[\dh(X,Y)] - \E[\dh(X',Y')] &= \P[X(e) \neq Y(e)] - \P[X'(e) \neq Y'(e)] \\
    &= \E \left[
        \begin{array}{c}
            \P\left[X(e) \neq Y(e) \middle| X(e'), Y(e') \text{ for } e' \neq e \right] \qquad \\
        - \; \P\left[X'(e) \neq Y'(e) \middle| X(e'), Y(e') \text{ for } e' \neq e \right]
        \end{array}\right].
\end{align}
Since, conditioned on $X(e'),Y(e')$ for all $e' \neq e$,
both $(X(e),Y(e))$ and $(X'(e),Y'(e))$ are pairs of Bernoullis with the same marginals
and $(X'(e),Y'(e))$ are always optimally coupled, the expression in the expectation on the right-hand
side is $\geq 0$.
Moreover, since we have assumed that $(X(e),Y(e))$ are \emph{not} optimally coupled for some realization
of $\{ X(e'),Y(e') \text{ for } e' \neq e \}$, the expectation is strictly positive.
This means $(X',Y')$ is a better coupling than $(X,Y)$, which is a contradiction to the assumption
that the optimal coupling $(X,Y)$ does not satisfy the conclusion of the lemma.
\end{proof}

%% file: sections/wasserstein/proof.tex
\subsection{Proof of main result}
\label{sec:wasserstein_proof}

With all of our preparation in the previous sections, we are finally ready to give a quick proof of
Theorem \ref{thm:wasserstein}, showing that
\begin{equation}
\label{eq:maingoal}
    \dW(\cM(n,p),\cG(n,p)) = \Theta(n^{3/2}).
\end{equation}

\begin{proof}[Proof of Theorem \ref{thm:wasserstein}]
First, note that again by \eqref{eq:wasdef_again} and the tower property,
under the optimal coupling between $X \sim \cM(n,p)$ and $Y \sim \cG(n,p)$
which minimizes the following expectation, we have
\begin{equation}
    \dW(\cM(n,p),\cG(n,p)) = \sum_{e \in \edgeset} \E \left[ \P \left[ X(e) \neq Y(e) \middle| X(e'), Y(e') \text{ for } e' \neq e \right] \right].
\end{equation}
Now by Lemma \ref{lem:optimality}, for any $e \in \edgeset$, under the optimal coupling, the distribution of $(X(e),Y(e))$
given all $X(e')$ and $Y(e')$ for $e' \neq e$ is that of two optimally coupled Bernoulli random variables
with means
\begin{equation}
    \E \left[ X(e) \middle| X(e') \text{ for } e' \neq e \right]
    \qquad \text{and} \qquad p
\end{equation}
respectively.
By the definition of $\cM(n,p)$, the expression on the left is equal to
\begin{equation}
    \phi(n^2 \partial_e \sH(X)) + \mathsf{error}(X),
\end{equation}
for some $\mathsf{error}(x) \lesssim 1$ which is zero unless $x$ is near the boundary of $\pball$.
In particular, $\mathsf{error}(x) = 0$ if $x \in \phalfball$, and thus $\E[\mathsf{error}(X)] = e^{-\Omega(n^2)}$.
In addition, since $p = \phi(2 \sH'(p))$ as discussed in Section \ref{sec:inputs_mixing}, we find that
\begin{equation}
    \dW(\cM(n,p),\cG(n,p)) = \sum_{e \in \edgeset} \E \left[
        \left|
            \phi(n^2 \partial_e \sH(X)) - \phi(2 \sH'(p))
        \right|
    \right] + e^{-\Omega(n^2)},
\end{equation}
and thus, since the summand above does not depend on $e$ by symmetry of $\cM(n,p)$, in order to show \eqref{eq:maingoal} it suffices to prove
the following bound for any fixed $e \in \edgeset$:
\begin{equation}
    \E \left[ \left| \phi(n^2 \partial_e \sH(X)) - \phi(2\sH'(p)) \right| \right] = \Theta(n^{-1/2}).
\end{equation}
Now let us Taylor expand $\phi$ to second order as follows:
\begin{equation}
\label{eq:taylorlast}
    \phi(n^2 \partial_e \sH(X)) - \phi(2\sH'(p)) = \phi'(2\sH'(p)) \cdot \left( n^2 \partial_e \sH(X) - 2 \sH'(p) \right)
    + O \left( \left( n^2 \partial_e \sH(X) - 2 \sH'(p) \right)^2 \right).
\end{equation}
For the second-order term,
By Corollary \ref{cor:concentration_marginal} with $k=2$ we have
\begin{equation}
    \E \left[ \left( n^2 \partial_e \sH(X) - 2 \sH'(p) \right)^2 \right]
    \lesssim n^{-1},
\end{equation}
and so taking the expectation of \eqref{eq:taylorlast} we obtain
\begin{equation}
    \E \left[ \left| \phi(n^2 \partial_e \sH(X)) - \phi(2 \sH'(p)) \right| \right]
    = \phi'(2\sH'(p)) \cdot \E \left[ \left| n^2 \partial_e \sH(X) - 2 \sH'(p) \right| \right]
    + O \left( n^{-1} \right).
\end{equation}
Thus, since $\phi'(2\sH'(p)) > 0$, it suffices to show that
\begin{equation}
\label{eq:innergoal}
    \E \left[ \left| n^2 \partial_e \sH(X) - 2 \sH'(p) \right| \right] = \Theta(n^{-1/2}).
\end{equation}
By Corollary \ref{cor:concentration_marginal}, the above expectation is $\lesssim n^{-1/2}$,
and on the other hand by Corollary \ref{cor:dhfluct}, it is $\gtrsim n^{-1/2}$.
This proves \eqref{eq:innergoal} and finishes the proof of Theorem \ref{thm:wasserstein}.
\end{proof}

%% file: sections/wasserstein/alternative.tex
\subsection{Alternative proof of lower bound via vertex degree CLT}
\label{sec:wasserstein_alternative}

In this section we give an alternative proof of the lower bound in Theorem \ref{thm:wasserstein},
via the quantitative central limit theorem for the degree of a vertex given by \cite[Theorem 1.3]{winstein2025quantitative}.
The idea for this proof was communicated to the author by Chris Jones following a discussion after a talk
given by the author at UC Davis, and we spell out the details here with his permission.

Let us first recall the statement of this central limit theorem; note that for the full version (not assuming the
phase-uniqueness-or-forest assumption as in \cite{winstein2025quantitative}), we refer to the discussion of
\cite[Section 5.2]{fkg}.
In the following, $\cV_j$ denotes the vertex set of the graph $G_j$ in the ERGM specification.
Additionally, for a vertex $\fu \in \cV_j$, $\sd_\fu$ denotes the degree of $\fu$ in $G_j$.
Finally, $\dWR$ denotes the $1$-Wasserstein distance between $\R$-valued random variables, i.e.\ the minimal expected absolute difference
under the optimal coupling.

\begin{theorem}[Theorem 1.3 of \cite{winstein2025quantitative}, Section 5.2 of \cite{fkg}]
\label{thm:degclt}
Define
\begin{equation}
    \varsigma_n^2 = \left( 1 - p(1-p) \sum_{j=1}^K \beta_j p^{\se_j-2} \sum_{\fu \in \cV_j} \sd_\fu (\sd_\fu-1) \right)^{-1} \times p (1-p) (n-1).
\end{equation}
Then for any deterministic vertex $v \in [n]$ and any $\eps > 0$, for $X \sim \cM(n,p)$ we have
\begin{equation}
    \dWR \left( \fr{\deg_v(X) - \E[\deg_v(X)]}{\varsigma_n}, Z \right) \lesssim n^{-\fr{1}{4}+\eps},
\end{equation}
where $Z \sim \cN(0,1)$.
\end{theorem}

Note that under the nondegeneracy assumption \eqref{eq:nondegen}, the variance proxy $\varsigma_n^2$ for $\deg_v(X)$ is larger (by a constant factor) than
the variance $p (1-p) (n-1)$ of the degree of a vertex in $\cG(n,p)$, since $\sd_\fu (\sd_\fu - 1)$ is nonzero for some vertex $\fu \in \cV_j$
if $G_j$ has any pair of adjacent edges.
Recall also that we may quantify the CLT for the vertex degree in $\cG(n,p)$ (which is a sum of $n-1$ i.i.d.\ $\Ber(p)$ variables)
by a standard application of Stein's method, e.g. \cite[Theorem 3.2]{ross2011fundamentals},
so that for $Y \sim \cG(n,p)$ and any $v \in [n]$ we have
\begin{equation}
\label{eq:gnpdegclt}
    \dWR \left( \fr{\deg_v(Y) - \E[\deg_v(Y)]}{\sqrt{p(1-p)(n-1)}}, Z' \right) \lesssim n^{-\fr{1}{2}},
\end{equation}
where again $Z' \sim \cN(0,1)$.

Now, by the hand-shaking lemma applied to the symmetric difference of the two graphs,
under any coupling between $X \sim \cM(n,p)$ and $Y \sim \cG(n,p)$, we have
\begin{align}
    \E \left[ \dh(X,Y) \right]
    &= \frac{1}{2} \sum_{v \in [n]} \E\left[ \sum_{e \ni v} |X(e) - Y(e)| \right] \\
    &\geq \frac{1}{2} \sum_{v \in [n]} \E \left[ |\deg_v(X) - \deg_v(Y)| \right],
\label{eq:hs}
\end{align}
applying the triangle inequality in the last step since $\deg_v(x) = \sum_{e \ni v} x(e)$.
Now by Theorem \ref{thm:marginal} we have
\begin{equation}
    \E[\deg_v(X)] = p(n-1) + O(1) = \E[\deg_v(Y)] + O(1),
\end{equation}
so, also using symmetry of the terms in \eqref{eq:hs}, we find that for any $v \in [n]$,
\begin{equation}
\label{eq:dhlast}
    \E \left[ \dh(X,Y) \right]
    \geq \frac{n}{2} \cdot \E\left[
        \left| \overline{\deg}_v(X) - \overline{\deg}_v(Y) \right|
    \right] - O(n),
\end{equation}
where we have set $\overline{\deg}_v(X) = \deg_v(X) - \E[\deg_v(X)]$ and similarly for $Y$.
Thus it remains to show that the expectation in the right-hand side above is $\gtrsim \sqrt{n}$ for any coupling between $X$ and $Y$.
For any coupling, this expectation is at least as large as $\dWR(\overline{\deg}_v(X), \overline{\deg}_v(Y))$ by the definition of the $1$-Wasserstein distance on $\R$.
Now by the triangle inequality for $\dWR$ we have
\begin{align}
    \dWR\left(
        \overline{\deg}_v(X),
        \overline{\deg}_v(Y)
    \right)
    &\geq
    \dWR\left(
        \varsigma_n Z, \sqrt{p(1-p)(n-1)} Z'
    \right) \\
    &\qquad -
    \dWR\left(
        \overline{\deg}_v(X), \varsigma_n Z
    \right) \\
    &\qquad -
    \dWR\left(
        \overline{\deg}_v(Y), \sqrt{p(1-p)(n-1)} Z'
    \right),
\end{align}
where, as above, $Z$ and $Z'$ are standard normal random variables.
Now by Theorem \ref{thm:degclt} and \eqref{eq:gnpdegclt} the second and third lines above are $O(n^{\frac{1}{4}})$ and $O(1)$ respectively.
On the other hand, the Wasserstein distance between two mean-zero Gaussians of variance $\sigma_1^2$ and $\sigma_2^2$ is exactly $| \sigma_1 - \sigma_2 | \sqrt{\frac{2}{\pi}}$,
as can be proven using the classical formula for $\dWR$ in terms of the inverse CDFs of the variables in question (for a derivation one may consult \cite[Corollary IV.3]{chhachhi20231}).
Thus, since
\begin{equation}
    \left| \varsigma_n - \sqrt{p(1-p)(n-1)} \right|
    \gtrsim \sqrt{n},
\end{equation}
we find that
\begin{equation}
    \dWR\left(
        \overline{\deg}_v(X),
        \overline{\deg}_v(X)
    \right) \geq c \sqrt{n} - C n^{\frac{1}{4}} - C
\end{equation}
for some constants $C, c > 0$.
Plugging this into \eqref{eq:dhlast} finishes the alternative proof of the lower bound.